\documentclass{article}

\usepackage{PRIMEarxiv}
\usepackage[utf8]{inputenc} 
\usepackage[T1]{fontenc}    
\usepackage{hyperref}       
\usepackage{url}            
\usepackage{booktabs}       
\usepackage{amsfonts}       
\usepackage{fancyhdr}       
\usepackage{graphicx}       
\graphicspath{ {./figures/} }

\usepackage{comment}
\usepackage{graphicx}
\usepackage{booktabs}
\usepackage{tikz}
\usetikzlibrary{positioning,arrows.meta,calc}
\usepackage{subcaption}
\usepackage{multirow}
\usepackage{algorithm}
\usepackage{algorithmic}

\usetikzlibrary{arrows.meta,positioning,fit,backgrounds,shapes.geometric}
\usepackage{natbib}

\RequirePackage{tgtermes}
\RequirePackage{newtxtext}
\RequirePackage{newtxmath}
\RequirePackage{bm}
\RequirePackage{endnotes}

\allowdisplaybreaks

\usepackage[nameinlink]{cleveref}
\crefformat{figure}{#2(figure~#1)#3}
\Crefformat{figure}{#2Figure~#1#3}
\crefformat{table}{#2(table~#1)#3}
\Crefformat{table}{#2Table~#1#3}
\crefformat{equation}{#2(~#1)#3}
\crefname{appendix}{appendix}{appendices}
\Crefname{appendix}{Appendix}{Appendices}
\crefformat{appendix}{#2appendix~#1#3}
\Crefformat{appendix}{#2Appendix~#1#3}

\title{Learning-Enhanced Reduce--Repair--Refine--Optimize Framework for Real-Time Train Platforming}

\author{
  Md Tabish Haque, Sebastian Sehmisch, Christopher Szymula, Nikola Be\v{s}inovi\'{c} \\
  Chair of Railway Operations, \\
  Technische Universit\"at Dresden \\
  Dresden, Germany\\[3pt]
  \texttt{\{mohammed\_tabish.haque, sebastian.sehmisch\}@tu-dresden.de} \\
  \texttt{\{christopher.szymula, nikola.besinovic\}@tu-dresden.de}
}

\usepackage{setspace}
\begin{document}
\maketitle

\begin{abstract}
Real-time train platforming under disruptions requires rapid revision of routing, platforming, and timing decisions. Dense traffic and shared infrastructure create tightly coupled conflicts, making real-time optimization computationally demanding. This paper proposes a learning-enhanced Reduce–Repair–Refine–Optimize framework that reduces and selectively relaxes the routing search space to balance computational efficiency, feasibility, and solution quality. The reduction stage constructs a restricted MILP by fixing selected routing decisions using route-specific predictive models and thresholds, thereby accounting for heterogeneity across routes. Two adaptive mechanisms subsequently modify the reduced search space: repair selectively releases route fixings to restore feasibility, while refinement releases additional fixings to recover routing alternatives and improve solution quality without reverting to the unrestricted problem. 
Computational experiments on real-life disrupted station-area instances show that the framework matches or improves upon a recent rTPP heuristic in more than 80\% of instances within 60~s. Under the common 180~s budget, the framework obtains strictly better terminal incumbents than the unrestricted MILP in 56.7\% of all instances and 77.4\% of hard instances. Moreover, it obtains better terminal incumbents than an otherwise identical variant without refinement in 71.9\% of hard instances. The results demonstrate that learning-based search-space reduction is most effective when the reduced search space remains adaptable during optimization, enabling high-quality feasible solutions within the real-time computation budget.
\end{abstract}%



\maketitle

%
\section{Introduction}\label{intro} 
Train station areas are particularly vulnerable to disruptions because dense traffic, complex track layouts, and limited routing and platform capacity create strong interactions among train movements. Disruptions can render the nominal timetable infeasible and require rapid revision of train routing, platforming, sequencing, and timing decisions. The resulting station-area rescheduling problem is commonly referred to as the real-time train platforming problem (rTPP) \citep{garcia2024efficient}. The objective is to obtain a feasible and conflict-free operating plan within the limited computational budget after a disruption, while limiting deviations from the nominal timetable, including delays, platform changes, route changes, and cancellations. The computational complexity of the rTPP arises from the strong interdependence among routing, platforming, and timing decisions. Each additional routing alternative expands the set of potential route-platform assignments and train interactions, which in turn constrains feasible timing decisions. Consequently, finding a high-quality operating plan within the limited computation time available constitutes a challenging combinatorial optimization problem (COP). 

Train routing and platforming have been studied across planning levels: strategic infrastructure capacity assessment \citep{zwaneveld1996routing}, tactical routing and platform assignment \cite{billionnet2003using, caprara2011solution, lusby2011routing}, and operational applications, including platform assignment in response to short-term demand fluctuations \citep{yang2026short} and real-time rescheduling under disruptions \citep{lusby2013set, garcia2024efficient}. Within the real-time setting, the need to obtain high-quality feasible solutions under tight computational budgets has motivated heuristic and metaheuristic methods, search-space reduction strategies, and emerging learning-based approaches for the rTPP \citep{garcia2024efficient, sehmisch2025data, sehmisch2026learning, zhang2026learning, gao2026train}. A recurring challenge is to obtain high-quality solutions within the limited computation time while retaining the detailed infrastructure representation and operational constraints needed to produce a feasible and conflict-free plan.
A growing body of research has explored the use of machine learning (ML) in COPs to guide the optimization process and reduce computational effort \citep{bengio2021machine, gasse2022machine, hizir2026large, rashedi2025machine}. Learning-based methods have shown promising results across several problem classes, including the traveling salesman problem, maximum clique enumeration, and fixed charge transportation problem \citep{sun2021generalization, tayebi2022learning, lauri2023learning, spieckermann2025reduce}. However, their effectiveness is strongly problem-dependent \citep{ding2020accelerating, gasse2019exact}, and these results do not directly transfer to other COPs, including the rTPP. 

For station-area rescheduling, ML-induced variable fixing has recently been explored to reduce the routing search space by predicting the stability of routes prescribed in the nominal timetable \citep{sehmisch2025data, sehmisch2026learning}. A route is considered \emph{stable} if all trains nominally assigned to it retain that route in the corresponding rescheduled solution. While these studies demonstrate the potential of learning-based search-space reduction, its use in the rTPP raises three important challenges: (i) route-specific variation in predictability, such that no single ML model performs uniformly well across routing decisions; (ii) infeasibility when prediction errors induce incorrect fixings that over-restrict the optimization problem; and (iii) solution-quality loss when the restricted problem remains feasible but excludes routing alternatives that could yield high-quality solutions.

Motivated by these challenges, we propose a multi-stage learning-and-optimization framework for the rTPP, referred to as the reduce--repair--refine--optimize (3RO) paradigm. The framework targets the routing decisions, which constitute a major source of combinatorial complexity. The reduction stage uses supervised ML at the route level to identify nominal routes expected to remain stable and fixes the associated routing binaries, thereby reducing the routing search space \citep{sama2016ant, nygren2023scope, sehmisch2025data}. Because ML-induced fixings may over-restrict the search space, the subsequent stages adapt this restriction through selective release of fixings. The repair stage releases a subset of route fixings to restore feasibility when the reduced model is infeasible, whereas the refinement stage releases a targeted subset of additional fixings to expand the restricted feasible region and recover routing alternatives that may improve solution quality, without reverting to the unrestricted mixed-integer linear program (MILP). The resulting restricted MILP is then optimized within the available real-time computation budget.

The main contributions of the paper are as follows: (1) We propose a learning-enhanced 3RO framework for the rTPP that integrates predictive route fixing, feasibility repair, refinement, and optimization within a common real-time computation budget. (2) We formulate route-stability prediction as a heterogeneous binary-relevance task in which the model class and decision threshold are selected separately for each route label, thereby accounting for route-specific heterogeneity in predictive performance. (3) We introduce selective release of ML-induced route fixings as an adaptive mechanism for addressing two distinct consequences of learning-based search-space reduction: repair restores feasibility when the restricted MILP is infeasible, whereas refinement releases additional fixings to recover routing alternatives for solution-quality improvement. (4) We characterize the downstream effects of learning-based reduction and evaluate the resulting computational trade-offs.  The analysis shows that infeasibility depends more strongly on which routing decisions are incorrectly fixed than on the aggregate number of prediction errors, with the fixings requiring release concentrating on a small subset of routes. We further compare alternative repair strategies, quantify the contribution of refinement to primal progress and terminal solution quality, and benchmark 3RO against both a recent rTPP heuristic and the unrestricted MILP.  Under the common 180~s budget, 3RO obtains strictly better incumbents than the unrestricted MILP in 56.7\% of all instances and 77.4\% of the hard instances.

\section{Literature Review}
This study relates to two main research streams: real-time train platforming and learning-enhanced solution methods for COPs. The first stream focuses on solution approaches for the rTPP and draws selective broader rRTMP literature for solution mechanisms, such as decomposition, rerouting search, and route selection, that are relevant to managing the computational complexity of station-area rescheduling despite differences in problem scope and decision structure (\Cref{sec:lit_rTPP}). The second stream reviews learning-enhanced solution methods for COPs and provides the methodological context for how ML can support computationally challenging optimization problems (\Cref{sec:lit_ML_4_COP}). Finally, \Cref{sec:lit_gaps} summarizes the resulting research gaps.

\subsection{Real-time train platforming and related solution approaches}\label{sec:lit_rTPP}
Early work on rTPP concentrated on platform reassignment and retiming, using exact optimization or heuristic search without explicitly optimizing train routes \citep{chakroborty2008optimum, zhang2020fast}. With routing incorporated, computational complexity has been mitigated through the decomposition of interrelated decisions or by limiting the optimization scope. \citet{ning2019rtpp} construct an iterative decomposition into a platform-route assignment with predetermined train schedules and later adjust for conflicting trains, whereas \citet{bai2021real} implement event-driven feasibility assessment followed by route and platform reassignment. Alternative approaches employ temporal decomposition or sequential construction. \citet{teng2024train} utilize a rolling-horizon approach that addresses smaller subproblems over successive time periods, while \citet{garcia2024efficient} develop a greedy constructive heuristic that sequentially optimizes train routing and scheduling while preserving earlier decisions.

In addition to rTPP-specific methodologies, associated rRTMP research decreases computational complexity by reducing the decision space examined during optimization. \citet{acuna2011mip} integrate variable fixing with local branching to constrain deviations from the nominal timetable. \citet{fischetti2017using} utilize heuristic bound tightening and variable fixing to refine variable domains and reduce the routing alternatives available to each train. At the scope level, \citet{van2018considering} establish a dynamic impact zone that focuses on conflict resolution on relevant trains and associated conflicts. Route-set reduction is explicitly outlined in the Train Routing Selection Problem (TRSP), which identifies a reduced set of alternate routes for each train prior to addressing the rRTMP \citep{sama2016ant, pascariu2022effective, pascariu2024formulation}. A related iterative method is presented in ROMA-based neighborhood search, which solves scheduling problems with predetermined routes and selectively modifies routing decisions through rerouting \citep{d2008reordering, corman2011effectiveness, sama2017variable}.

Learning-based methods shift some computational effort to offline training and use the learned information to accelerate online rescheduling. In the context of rTPP, \citet{zhang2026learning} employ a deep reinforcement learning policy to generate platforming and rescheduling decisions. Related studies in railway rescheduling combine learning with optimization by fixing selected decisions prior to optimizing, using reinforcement learning and supervised learning, respectively \citep{zhang2024integrated, su2024data}. \citet{nygren2023scope} define a \emph{core problem} at a high level of aggregation, where the optimization is limited to trains and related decisions that need to be rescheduled, while decisions outside the scope are fixed. More directly for the rTPP, \citet{sehmisch2025data, sehmisch2026learning} implement a learning-driven approach to route reduction by predicting the stability of nominal routes in the rescheduled solution and subsequently adjusting the relevant route-selection binaries prior to optimization.

Overall, existing studies improve computational tractability through decomposition, sequential solution construction, and restriction of the rescheduling problem or its search space. The effectiveness of these methods depends on how much coordination and decision flexibility are preserved. Decomposition requires coordination across subproblems, sequential construction may limit the revision of prior decisions, and variable fixing may exclude alternatives needed for feasibility or objective improvement. The common challenge is therefore to reduce computational effort while preserving sufficient coordination and flexibility to obtain high-quality feasible solutions within real-time computation budgets.

\subsection{Learning-enhanced solution methods for COPs}\label{sec:lit_ML_4_COP}
Learning-enhanced optimization integrates ML into an explicit optimization model or solution algorithm, either to support internal solver decisions or to restrict the optimization problem by identifying promising decision components \citep{bengio2021machine, scavuzzo2024machine, sehmisch2025data}. The latter, often referred to as learning-based problem reduction, is particularly relevant for time-critical combinatorial optimization: predictions from previously solved instances are used to remove, fix, or prioritize decision components before the resulting restricted problem is solved.

Beyond the railway domain, learning-based reduction has been applied by predicting relevant decision components and restricting the subsequent optimization accordingly, including vertices in clique problems \citep{sun2019using,lauri2023learning}, edges in routing problems \citep{fitzpatrick2023learning}, and partial variable assignments in MILPs \citep{ding2020accelerating,nair2020solving}. More recently, \citet{spieckermann2025reduce} formalize a reduce--then--optimize approach for the fixed-charge transportation problem and augment the predicted edge set with a greedily constructed feasible solution. Although this construction preserves feasibility before optimization, their results show that strong predictive scores can coexist with large optimality gaps when a solution edge falls just outside the retained set. In airline disruption recovery, \citet{rashedi2025machine} prune the aircraft-recovery search space, while \citet{hizir2026large} impose learned crew-recovery restrictions through soft cuts, as hard enforcement was found to produce infeasible or low-quality solutions. The latter further show that downstream solution quality is more closely related to correctly identifying a small subset of solution-relevant decisions than to aggregate predictive performance. Related studies have also introduced explicit mechanisms to address infeasible learned decisions, including infeasibility elimination before variable fixing \citep{yilmaz2024expandable} and optimization-based correction of ML predictions \citep{ojha2025outbound}.

Overall, these studies show that learning-based reduction can improve computational tractability, but aggregate predictive metrics alone do not fully characterize its effects on the subsequent optimization. These effects depend not only on the number of prediction errors but also on which decisions they affect.

\subsection{Scientific gaps}\label{sec:lit_gaps}
The reviewed literature highlights several gaps that motivate the present study and are as follows: 
\begin{enumerate}
\item Existing rTPP and related railway rescheduling methods reduce computational effort through decomposition, sequential solution construction, and restrictions on the optimization scope or routing alternatives. Learning-based reduction shifts part of the restriction-selection process to offline prediction. However, limited attention has been given to adapting the resulting reduced search space during optimization when the imposed restrictions render the problem infeasible or exclude routing alternatives that could improve solution quality.

\item Evaluations of learning-based reduction typically emphasize aggregate predictive metrics, runtime, and terminal solution quality. Evidence from related COPs nevertheless shows that similar predictive performance can lead to different optimization outcomes depending on which decisions are predicted incorrectly. Consequently, how prediction errors induce infeasibility, which fixings must be released to restore feasibility, and how recovering routing alternatives affects primal progress and terminal solution quality under a real-time computation budget remain insufficiently characterized in railway rescheduling.

\item Learning-based reduction often applies a common model class and decision threshold across decision components. Although existing approaches recognize differences in predictability through component selection or classifier aggregation, component-specific selection of both the predictive model and decision threshold remains insufficiently examined. This is particularly relevant for the rTPP, where route labels may differ substantially in class balance and predictive performance.
\end{enumerate}

\section{Problem Description}\label{problemdescription} 
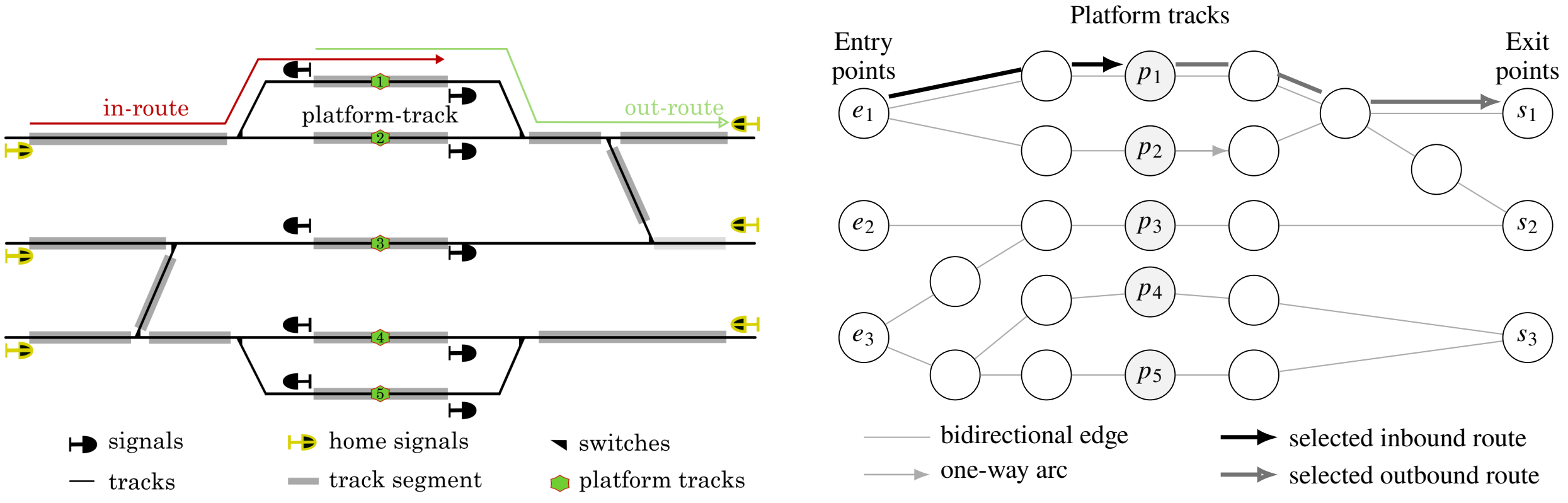
\begin{figure}[b]
    \centering
    \begin{subfigure}[t]{0.48\textwidth}
        \centering
        \resizebox{\linewidth}{!}{%
        \includegraphics{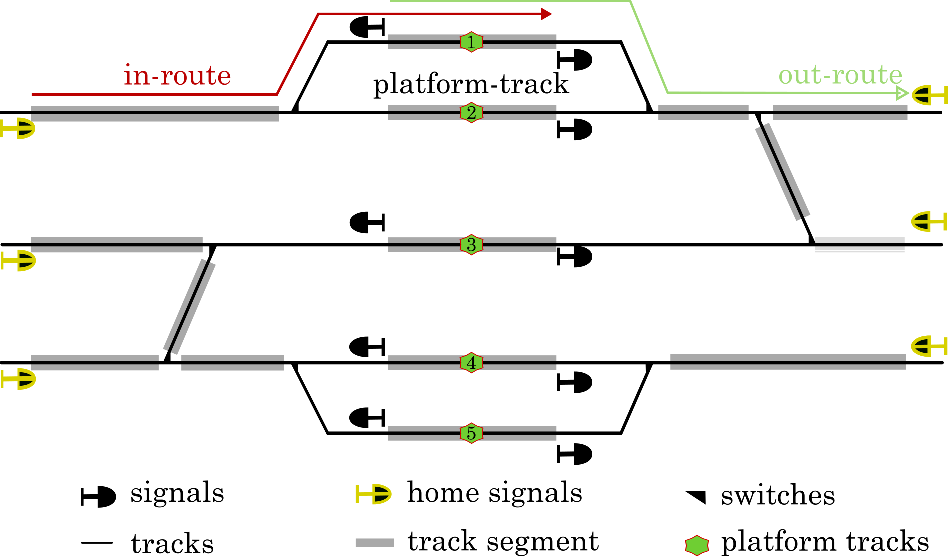}}
        \caption{A station area representation}
        \label{fig:infra_model}
    \end{subfigure}
    \hfill
    \begin{subfigure}[t]{0.48\textwidth}
    \centering
    \resizebox{\linewidth}{!}{%
\begin{tikzpicture}[
    >=Latex,
    scale=1.0,
    every node/.style={font=\small},
    track/.style={circle, draw=black, fill=white, minimum size=6mm, inner sep=0pt},
    platform/.style={circle, draw=black, fill=gray!10, minimum size=6mm, inner sep=1pt},
    boundary/.style={circle, draw=black, fill=white, minimum size=6mm, inner sep=1pt},
    infra/.style={line width=0.45pt, draw=gray!65},
    infraarc/.style={->, line width=0.55pt, draw=gray!65},
    inroute/.style={line width=1.45pt, draw=black,
                    -{Latex[length=2.6mm, width=2.0mm]}},
    outroute/.style={line width=1.45pt, draw=black!55,
                     -{Latex[length=2.6mm, width=2.0mm, open, fill=white]}},
    routeshift/.style={transform canvas={yshift=4pt}},
    move/.style={->, line width=0.9pt, draw=black},
    labeltext/.style={font=\small},
    legendtext/.style={font=\footnotesize, anchor=west}
]

\node[boundary] (entry1) at (0,1.35) {$e_1$};
\node[boundary] (entry2) at (0,0) {$e_2$};
\node[boundary] (entry3) at (0,-1.35) {$e_3$};

\node[boundary] (exit1) at (8.0,1.35) {$s_1$};
\node[boundary] (exit2) at (8.0,0.0) {$s_2$};
\node[boundary] (exit3) at (8.0,-1.35) {$s_3$};

\node[track] (a3) at (1.1,-0.675) {};
\node[track] (a4) at (1.1,-1.8) {};

\node[track] (b1) at (2.2,1.8) {};
\node[track] (b2) at (2.2,0.9) {};
\node[track] (b3) at (2.2,-0) {};
\node[track] (b4) at (2.2,-0.9) {};
\node[track] (b5) at (2.2,-1.8) {};

\node[platform] (p1) at (3.45,1.8) {$p_1$};
\node[platform] (p2) at (3.45,0.9) {$p_2$};
\node[platform] (p3) at (3.45,0) {$p_3$};
\node[platform] (p4) at (3.45,-0.8) {$p_4$};
\node[platform] (p5) at (3.45,-1.8) {$p_5$};

\node[track] (c1) at (4.7,1.8) {};
\node[track] (c2) at (4.7,0.9) {};
\node[track] (c3) at (4.7,-0.0) {};
\node[track] (c4) at (4.7,-0.9) {};
\node[track] (c5) at (4.7,-1.8) {};

\node[track] (d1) at (5.8,1.35) {};

\node[track] (e1) at (6.9,0.675) {};

\draw[infra] (entry1) -- (b1);
\draw[infra] (entry1) -- (b2);
\draw[infra] (entry2) -- (b3);
\draw[infra] (entry3) -- (a3);
\draw[infra] (entry3) -- (a4);

\draw[infra] (a3) -- (b3);
\draw[infra] (a4) -- (b4);
\draw[infra] (a4) -- (b5);

\draw[infra] (b1) -- (p1);
\draw[infra] (b2) -- (p2);
\draw[infra] (b3) -- (p3);
\draw[infra] (b4) -- (p4);
\draw[infra] (b5) -- (p5);

\draw[infra] (p1) -- (c1);
\draw[infra] (p2) -- (c2);
\draw[infra] (p3) -- (c3);
\draw[infra] (p4) -- (c4);
\draw[infra] (p5) -- (c5);

\draw[infra] (c1) -- (d1);
\draw[infra] (c2) -- (d1);

\draw[infra] (d1) -- (e1);

\draw[infra] (d1) -- (exit1);
\draw[infra] (c3) -- (exit2);
\draw[infra] (e1) -- (exit2);
\draw[infra] (c4) -- (exit3);
\draw[infra] (c5) -- (exit3);

\draw[infraarc] (p2) -- (c2);


\draw[inroute, routeshift]
    (entry1) -- (b1) -- (p1);

\draw[outroute, routeshift]
    (p1) -- (c1) -- (d1) -- (exit1);

\node[label={[align=center]above:Entry\\points}] at (0,1.45) {};
\node[label={[align=center]above:Exit\\points}] at (8.0,1.45) {};
\node[label={[align=center]below:Platform tracks}] at (3.45,2.90) {};

\begin{scope}[shift={(0,-2.55)}]
    \draw[infra] (0,0) -- (0.8,0);
    \node[legendtext] at (0.8,0) {bidirectional edge};

    \draw[inroute] (4.3,0) -- (5.0,0);
    \node[legendtext] at (5.0,0){selected inbound route};

    \draw[infraarc] (0,-0.45) -- (0.8,-0.45);
    \node[legendtext] at (0.8,-0.45){one-way arc};

    \draw[outroute] (4.3,-0.45) -- (5.0,-0.45);
    \node[legendtext] at (5.0,-0.45) {selected outbound route};
\end{scope}

\end{tikzpicture}}
    \caption{Infrastructure-graph representation}
    \label{fig:infra_graph}
\end{subfigure}
    \caption{Illustrative infrastructure and graph-based representation of a station-area train movement}
\end{figure}

We consider a station area composed of track segments, switches, platform tracks, and signaling, as illustrated in \Cref{fig:infra_model}. The physical station area infrastructure is represented as a mixed graph $\mathcal{G}^{\mathcal{I}} = (\mathcal{V}, \mathcal{E}, \mathcal{A})$, where vertices $\mathcal{V}$ denote track segments, including platform tracks; undirected edges $\mathcal{E}$ represent bidirectional movements; and directed arcs $\mathcal{A}$ represent one-way movements between adjacent segments, as determined by the track layout, switch positions, and signaling rules (see \Cref{fig:infra_graph}). Let $\mathcal{V}^{\mathrm{pf}}\subseteq \mathcal{V}$ denote the platform tracks $p_k$, and let $\mathcal{V}^{\mathrm{b}}\subseteq \mathcal{V}$ the entry vertices $e_h$ and exit vertices $s_o$; the remaining vertices represent intermediate track segments. A feasible train path on $\mathcal{G}^{\mathcal{I}}$ is an ordered sequence of adjacent vertices connected by edges or arcs from an entry vertex $e_h$ to an exit vertex $s_o$, passing through a platform track $p_k$. Because routing and platforming decisions are jointly determined, we decompose each feasible train path at the platform track $p_k$ into an inbound and an outbound route. An \emph{inbound route} $r\in\mathcal{R}^{+}$ is a feasible path from an entry vertex to a platform track, and an \emph{outbound route} $r\in\mathcal{R}^{-}$ is a feasible path from a platform track to an exit vertex, with $\mathcal{R}=\mathcal{R}^{+}\cup\mathcal{R}^{-}$. We denote by $\mathcal{V}(r)\subseteq \mathcal{V}$ the set of vertices traversed by route $r$. An inbound and an outbound route are \emph{compatible} if they serve the same platform track. The resulting compatible inbound--outbound route pair defines a feasible train path from an entry vertex to an exit vertex via the common platform track.

The nominal timetable specifies, for each train $t\in \mathcal{T}$, its scheduled arrival and departure times, assigned platform track, and inbound-outbound route pair. For each route $r\in\mathcal{R}$, let $\mathcal{T}_r = \{t\in\mathcal{T} : r \text{ belongs to the inbound--outbound route pair assigned to } t \text{ in the nominal timetable}\}$. The set of nominal routes is then defined as $\mathcal{R}^{\mathrm{nom}} = \{r\in\mathcal{R} : \mathcal{T}_r\neq\emptyset\}$. For rescheduling, each train $t\in\mathcal{T}$ is associated with feasible inbound and outbound route sets $\mathcal{R}^{+}_t\subseteq\mathcal{R}^{+}$ and $\mathcal{R}^{-}_t\subseteq\mathcal{R}^{-}$, respectively, determined by its entry and exit vertices, direction, and platform-use restrictions. The compatible inbound--outbound route pairs in $\mathcal{R}^{+}_t\times\mathcal{R}^{-}_t$ constitute the routing and platforming alternatives available to train $t$. Routes in $\mathcal{R}\setminus\mathcal{R}^{\mathrm{nom}}$ may serve as rerouting alternatives during rescheduling. 

Disruptions are modeled as the temporary unavailability of vertices, edges, or arcs in $\mathcal{G}^{\mathcal{I}}$, representing blocked tracks and restricted movements caused by switch or signaling failures. Such disruptions reduce effective station capacity and may render the nominal timetable infeasible. The resulting rTPP determines revised routing, platforming, sequencing, and timing decisions to obtain a feasible and conflict-free operating plan under the reduced station capacity. Let $u=(\alpha,\beta)$ denote the routing binaries, where $\alpha^r_t=1$ if train $t$ is assigned inbound route $r\in\mathcal{R}^{+}_t$, and $\beta^r_t=1$ if train $t$ is assigned outbound route $r\in\mathcal{R}^{-}_t$. All remaining variables are collected in $v=(d,q)$, where $d$ contains operational-deviation variables representing cancellations, modeled by assigning a dummy route, as well as delays, platform changes, and route changes, while $q$ contains auxiliary timing and sequencing variables. To keep the presentation compact, we express the rTPP in block-matrix form. For each constraint family indexed by $\kappa$, $A_x^\kappa$ denotes the coefficient matrix associated with variable block $x$ and $b^\kappa$, where applicable, denotes the corresponding right-hand-side vector. The rTPP can then be formulated as
\begingroup
\setlength{\abovedisplayskip}{4pt}
\setlength{\belowdisplayskip}{4pt}
\setlength{\belowdisplayshortskip}{4pt}
\begin{align}
\text{(rTPP)} \quad 
\min_{u,v} \quad 
& c_d^{\top}d \label{eq:trp_obj}\\
\text{s.t.} \quad 
& A^{\mathrm{time}}_u u + A^{\mathrm{time}}_v v \geq b^{\mathrm{time}} \label{eq:trp_time}\\
& A^{\mathrm{flow}}_u u + A^{\mathrm{flow}}_d d = b^{\mathrm{flow}} \label{eq:trp_flow}\\
& A^{\mathrm{head}}_u u + A^{\mathrm{head}}_v v \geq b^{\mathrm{head}} \label{eq:trp_head}\\
& A^{\mathrm{dev}}_u u + A^{\mathrm{dev}}_v v \leq d \label{eq:trp_dev}\\
& u_i \in \{0,1\} \quad i \in B^r \label{eq:trp_routebin}\\
& v_i \in \{0,1\} \quad i \in B^o \label{eq:trp_otherbin}\\
& v_i \in \mathbb{Z}_{+} \quad i \in I \label{eq:trp_int}\\
& v_i \in \mathbb{R}_{+} \quad i \in C \label{eq:trp_cont}
\end{align}
\endgroup
The objective \eqref{eq:trp_obj} minimizes the weighted cost of operational deviations. The timing constraints \eqref{eq:trp_time} determine arrival and departure times within the station area and the resulting delays relative to the nominal timetable. The flow conservation constraints \eqref{eq:trp_flow} ensure that each train is either assigned a compatible pair of inbound and outbound routes or cancelled. The headway constraints \eqref{eq:trp_head} enforce infrastructure-capacity and conflict-avoidance requirements for inbound, outbound, and crossing movements. Finally, constraints \eqref{eq:trp_dev} define the deviation variables for platform and route changes relative to the nominal timetable. The index set $B^r$ corresponds to the routing binaries, whereas $B^o$, $I$, and $C$ define the binary, integer, and continuous components of $v$. The detailed MILP formulation of the station-area train rescheduling problem is provided in \hyperref[app:rtpp_formulation]{Appendix}. 

The computational cost of the rTPP is primarily driven by the number of routing variables and operational constraints \citep{pascariu2022effective}, together with the strongly coupled structure induced by shared infrastructure and train interactions. Additional inbound and outbound route alternatives enlarge the binary routing space, while higher traffic density increases the number of sequencing and headway interactions \citep{corman2011effectiveness}. These effects are amplified under disruptions, which reduce the available station capacity and increase competition for the remaining routes and platforms. However, disruptions are often spatially localized \citep{van2018considering, sehmisch2025data}, their operational effects are therefore typically concentrated on a subset of trains and route assignments. Consequently, many trains can retain their nominal routes, whereas rerouting is required primarily for trains affected directly or through propagated route, platform, and headway interactions. This structure provides the basis for route-stability prediction: nominal routes predicted to remain stable can be fixed, thereby reducing the routing search space while leaving the remaining routing decisions and all timing, sequencing, and deviation variables available for optimization.

\section{Methodology}\label{sec:methodology}
\begin{figure}[t]
    \centering
    \resizebox{0.95\linewidth}{!}{\resizebox{0.98\linewidth}{!}{%
\begin{tikzpicture}[
    font=\normalsize,
    >=Latex,
    box/.style={
        draw,
        rounded corners=2pt,
        align=center,
        minimum height=1.00cm,
        minimum width=2.25cm,
        inner sep=5pt
    },
    widebox/.style={
        draw,
        rounded corners=2pt,
        align=center,
        minimum height=1.00cm,
        minimum width=2.55cm,
        inner sep=5pt
    },
    decision/.style={
        draw,
        diamond,
        aspect=2.15,
        align=center,
        minimum width=2.25cm,
        minimum height=1.15cm,
        inner sep=2pt
    },
    stagebox/.style={
        draw,
        rounded corners=8pt,
        inner sep=10pt
    },
    arrow/.style={
        ->,
        thick,
        shorten >=1pt,
        shorten <=1pt
    },
    dasharrow/.style={
        ->,
        thick,
        dashed,
        shorten >=1pt,
        shorten <=1pt
    }
]

\coordinate (commoncenter) at (0,0);

\coordinate (offmid) at (0,0);

\node[widebox, left=5mm of offmid] (offsol)
    {Reference rTPP\\solutions};

\node[widebox, right=15mm of offmid] (labels)
    {Route-stability\\labels};

\node[widebox, left=15mm of offsol] (offdata)
    {Reference disruption\\instances};

\node[widebox, right=25mm of labels] (scoper)
    {Trained\\Scoper};

\node[
    font=\bfseries\large,
    above=9mm of $(offsol.north)!0.5!(labels.north)$
] (offtitle)
    {Offline learning};

\draw[arrow] (offdata.east) -- (offsol.west);
\draw[arrow] (offsol.east) -- (labels.west);
\draw[arrow] (labels.east) -- (scoper.west);

\node[widebox, below=43mm of offmid] (red)
    {Reduced MILP\\$\mathrm{rTPP}_{\mathrm{red}}$};

\node[widebox, left=7mm of red] (pred)
    {Route-stability\\prediction};

\node[widebox, left=7mm of pred] (inst)
    {Disrupted rTPP\\instance $j'$};

\node[decision, right=7mm of red] (feas)
    {$\mathrm{rTPP}_{\mathrm{red}}$\\feasible?};

\node[box, right=9mm of feas] (repair)
    {Auxiliary\\repair MILP};

\node[box, below=18mm of feas] (inc)
    {Feasible\\incumbent};

\node[widebox, left=7mm of inc] (refine)
    {Interaction-guided\\refinement};

\node[box, left=7mm of refine] (opt)
    {Optimize\\$\mathrm{rTPP}_{\mathrm{ref}}$};

\node[widebox, left=7mm of opt] (out)
    {Final revised\\routes and schedule};

\node[
    font=\bfseries\large,
    above=10mm of red.north
] (ontitle)
    {Online solution process};

\draw[arrow] (inst.east) -- (pred.west);
\draw[arrow] (pred.east) -- (red.west);
\draw[arrow] (red.east) -- (feas.west);

\draw[arrow]
    (feas.south)
    --
    node[right=2pt] {Yes}
    (inc.north);

\draw[arrow]
    (feas.east)
    --
    node[above=2pt] {No}
    (repair.west);

\draw[arrow]
    (repair.south)
    |-
    (inc.east);

\draw[arrow] (inc.west) -- (refine.east);
\draw[arrow] (refine.west) -- (opt.east);
\draw[arrow] (opt.west) -- (out.east);

\draw[arrow]
    (opt.south)
    -- ++(0,-7mm)
    -|
    (refine.south);
    
\draw[dasharrow]
    (scoper.south)
    -- ++(0,-12mm)
    -|
    (pred.north);

\coordinate (offlineLeft)  at (-8.25,0);
\coordinate (offlineRight) at ( 8.25,0);

\coordinate (onlineLeft)  at (-8.25,-5.4);
\coordinate (onlineRight) at ( 8.25,-5.4);
\coordinate (onlineBottom) at ($(opt.south)!0.5!(refine.south)+(0,-10mm)$);
\begin{scope}[on background layer]

\node[
    stagebox,
    fit={
        (offtitle)
        (offdata)
        (offsol)
        (labels)
        (scoper)
        (offlineLeft)
        (offlineRight)
    }
] (offlinebox) {};

\node[
    stagebox,
    fit={
        (ontitle)
        (inst)
        (pred)
        (red)
        (feas)
        (repair)
        (inc)
        (refine)
        (opt)
        (out)
        (onlineLeft)
        (onlineRight)
        (onlineBottom)     
    }
] (onlinebox) {};

\end{scope}

\end{tikzpicture}%
}}
    \caption{Overview of the 3RO framework, consisting an offline learning and online solution process for rTPP}
    \label{fig:offline_online}
\end{figure}
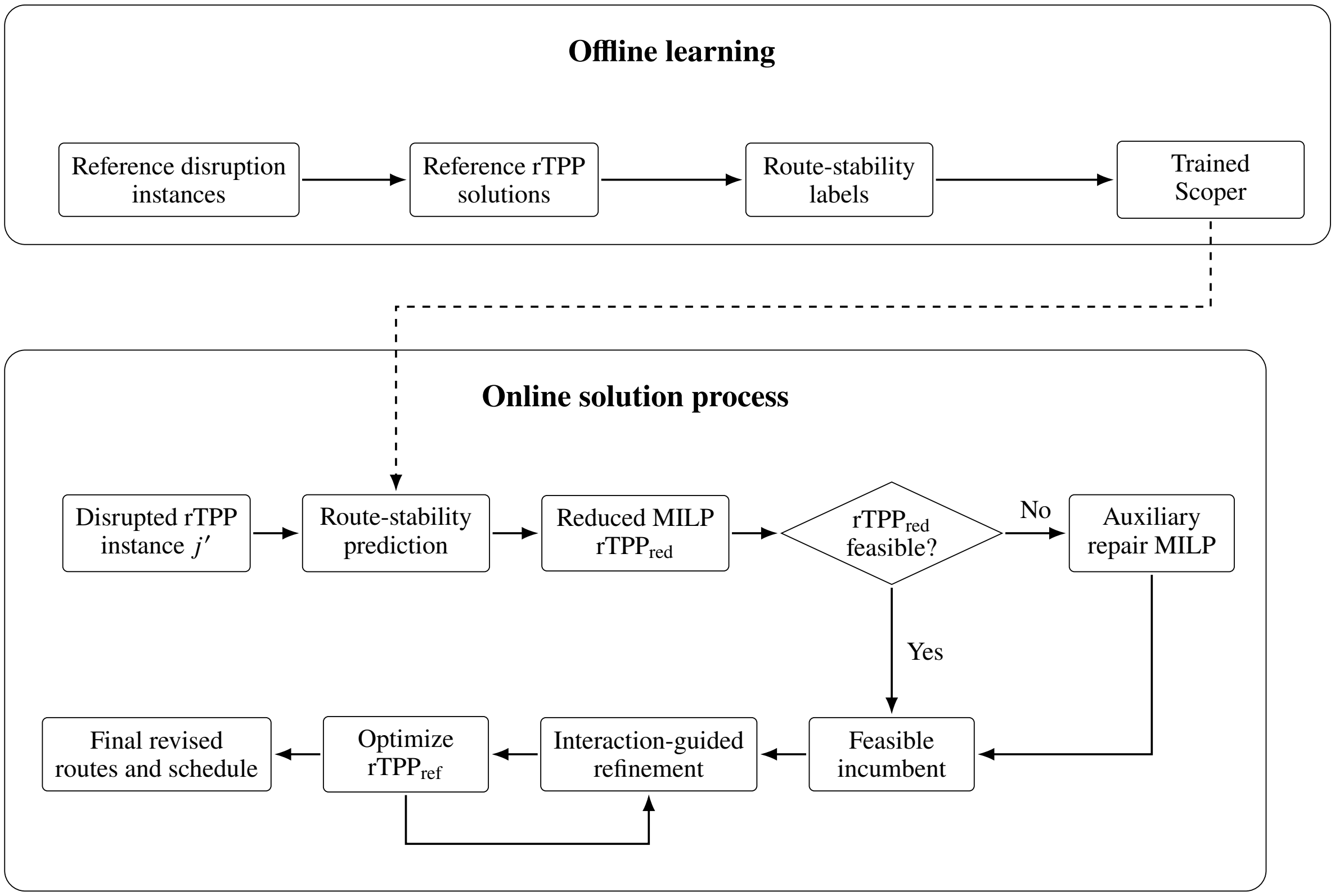

We propose a reduce--repair--refine--optimize (3RO) framework for solving the rTPP in real time. The framework combines learning-based search-space reduction with adaptive release of selected variable fixings to balance computational tractability, feasibility, and solution quality. It adopts an offline--online architecture, with model training conducted offline and reduction, repair, refinement, and optimization performed online. \Cref{fig:offline_online} provides an overview of this architecture and the information flow between its offline and online phases.

In the offline phase, reference disruption instances are solved to obtain rescheduling solutions from which route-stability labels are derived. These labels are used to train route-specific binary classifiers that collectively constitute the \emph{Scoper}. Online, the trained Scoper predicts which nominal routes remain stable for a disrupted rTPP instance. The reduction stage maps these predictions to a subset of routing binaries that are fixed to their nominal values, thereby reducing the routing search space and yielding the restricted MILP $\mathrm{rTPP}_{\mathrm{red}}$, with a smaller routing search space than the original rTPP (\Cref{sec:reduction}). Since ML-induced route fixings may render $\mathrm{rTPP}_{\mathrm{red}}$ infeasible, the repair stage is invoked when feasibility is not preserved (\Cref{sec:repair}). It solves an auxiliary repair problem to identify and release a limited subset of route fixings, thereby restoring feasibility while preserving most of the ML-induced reduction. Once feasibility is established, the resulting restricted model and its incumbent are passed to the refinement stage (\Cref{sec:refine}).  The remaining route fixings may still exclude routing alternatives that could better accommodate disruption-specific routing, platforming, and headway interactions and thereby improve solution quality. Refinement therefore selectively releases additional active route fixings to recover such alternatives, yielding the refined model $\mathrm{rTPP}_{\mathrm{ref}}$. This expands the reduced routing search space while the model remains restricted relative to the original rTPP. Finally, the optimization stage solves $\mathrm{rTPP}_{\mathrm{ref}}$ within the remaining real-time computation budget to obtain a high-quality feasible solution.

\subsection{Learning-based search-space reduction}\label{sec:reduction}
The reduction stage decreases the combinatorial complexity of the rTPP by fixing selected routing binaries $u=(\alpha,\beta)$. These variables are natural candidates for learning-based reduction because they define the discrete route-choice structure and largely determine the number of route-platform combinations and train interactions that must be resolved through headway and sequencing constraints. Route alternatives are indexed over the base infrastructure graph $\mathcal{G}^{\mathcal{I}}$, which provides a consistent route-label space across different disruption instances. The remaining variables $v=(d,q)$, including timing, sequencing, and operational-deviation decisions, remain endogenous MILP variables whose values are determined by the realized route, platform, and headway interactions. 

The resulting search-space reduction is determined by the Scoper's route-stability predictions. For each disruption instance, a binary outcome is predicted for every nominal route; the resulting prediction task is therefore formulated as a multi-label classification problem. We address this task using the binary-relevance strategy \citep{herrera2016multilabel}, under which a separate binary classifier is trained for each route label.

\subsubsection{Offline training}
The Scoper is trained offline using reference disruption instances and their corresponding rescheduling solutions obtained from the original rTPP, \eqref{eq:trp_obj}\text{–-}\eqref{eq:trp_cont}. Each instance is described by spatial, temporal, and operational features, including the disruption location, affected routes, topological characteristics of the disrupted element, disruption start time, and number of affected trains. These features are constructed from the disruption declaration and the nominal timetable and are therefore available when 3RO is invoked online. The corresponding rescheduling solutions provide the route-stability labels. For each disruption instance $j$ and nominal route $r\in \mathcal{R}^{\mathrm{nom}}$, the binary label $y_{jr}$ indicates whether all trains nominally assigned to $r$ retain that route in the reference solution for instance $j$:

\begin{equation*}
y_{jr} =
\begin{cases}
1, & \text{if route } r \text{ is retained,} \\
0, & \text{otherwise}.
\end{cases}
\end{equation*}

A positive label $y_{jr} = 1$, therefore, indicates that route $r$ remains stable under the disruption instance $j$, whereas a negative label $y_{jr} = 0$ indicates that at least one train nominally assigned to route $r$ requires rerouting. Since the rTPP may admit multiple optimal solutions, the labels are defined relative to the selected reference solution for each instance and may differ across alternative reference solutions. The resulting labels therefore provide a consistent, but not necessarily unique, representation of route stability.

Let $\mathcal{M}$ denote the set of candidate model classes. A separate classifier is trained for each model class $m\in\mathcal{M}$ and route label $r$ using five-fold cross-validation on the training data. Average precision (AP) is used as the performance criterion throughout offline training because it summarizes the precision--recall curve over the predicted probabilities without requiring a fixed decision threshold, which has not yet been selected at this stage. For a set of predictions ranked by decreasing predicted probability, AP is computed as $\mathrm{AP} = \sum_{n} (R_n - R_{n-1})\, P_n$, where $P_n$ and $R_n$ denote the precision and recall obtained when the $n$-th highest predicted probability is used as the cutoff, i.e., aggregating performance across all possible cutoffs rather than a single fixed decision threshold. Within the five-fold cross-validation procedure, AP is computed on the held-out predictions of each fold and averaged across folds to obtain the criterion used for recursive feature elimination and hyperparameter tuning. To reduce input dimensionality and improve generalization, recursive feature elimination \citep{guyon2003introduction, haider2024subnetwork} selects the most informative subset of the engineered features, while hyperparameters are tuned by grid search. For each trained classifier, two candidate decision thresholds are considered for converting predicted route-stability probabilities into binary predictions: the default threshold $\phi^{\mathrm{def}}$ and a route- and model-specific threshold $\phi_{r,m}^{\mathrm{ROC}}$, obtained from receiver operating characteristic (ROC) analysis \citep{fawcett2006introduction} on the cross-validated predictions. The ROC-based threshold balances the true-positive rate, defined as the proportion of stable cases correctly classified, against the false-positive rate, defined as the proportion of potentially affected cases incorrectly classified as stable. Together, these thresholds form the candidate set $\Phi_{r,m} = \{\phi^{\mathrm{def}}, \phi_{r,m}^{\mathrm{ROC}}\}$ used in the subsequent label-wise threshold selection procedure.

\subsubsection{Heterogeneous Scoper ensemble}
We construct the heterogeneous Scoper by selecting the model class and decision threshold separately for each route label, thereby accommodating variation in predictive performance across labels. Model and threshold selection follows a two-step procedure using a validation set disjoint from the training data. First, the model class attaining the highest label-specific average precision ($\mathrm{AP}_r$) is selected: 
\begin{equation}
m^*{(r)}=\arg\max_{m\in\mathcal{M}} \mathrm{AP}_{r}(m),
\end{equation}

Here, $\mathrm{AP}_{r}(m)$ measures the threshold-independent ranking performance of model class $m$ for route label $r$; the same criterion is used for feature selection and hyperparameter tuning. Conditional on the selected model class $m^{*}(r)$, the threshold-selection step determines 
\begin{equation}
\phi^*(r)=\arg\max_{\phi\in\Phi_{r,m^*(r)}} F_{1,r}\bigl(m^*(r),\phi\bigr),
\end{equation}
where $F_{1,r}$ denotes the $F_1$-score for route label $r$, computed as the harmonic mean of precision and recall \citep{herrera2016multilabel}. Ties in the model-selection step are resolved by applying the threshold-selection criterion to each tied model class and retaining the resulting model--threshold pairs with the highest $F_1$-score. If multiple model classes attain the same $F_1$-score, the tie is resolved using a predefined model-class priority, discussed in \Cref{sec:experimental_design}. The resulting heterogeneous Scoper is $\mathcal{S}=\{f_r^{m^*(r),\phi^*{(r)}}: r\in\mathcal{R}^{\mathrm{nom}}\}$, where $f_r^{m^*(r),\phi^*(r)}$ denotes the route-specific classifier and its associated decision threshold. The classifiers and thresholds are selected offline and remain fixed during online application.

\subsubsection{Online route fixing and restricted rTPP}
In the online phase, the trained heterogeneous Scoper is applied to a newly observed disruption instance $j'$. For each route label $r$, the route-specific classifier produces a stability probability $\hat{p}_{j'r}\in[0,1]$, which is converted into a binary route-stability prediction using the corresponding decision threshold $\phi^*{(r)}$: 
\begin{equation*}
\hat{y}_{j'r} =
\begin{cases}
1, & \text{if } \hat{p}_{j'r}\geq \phi^{*}(r),\\
0, & \text{otherwise}.
\end{cases}
\end{equation*}

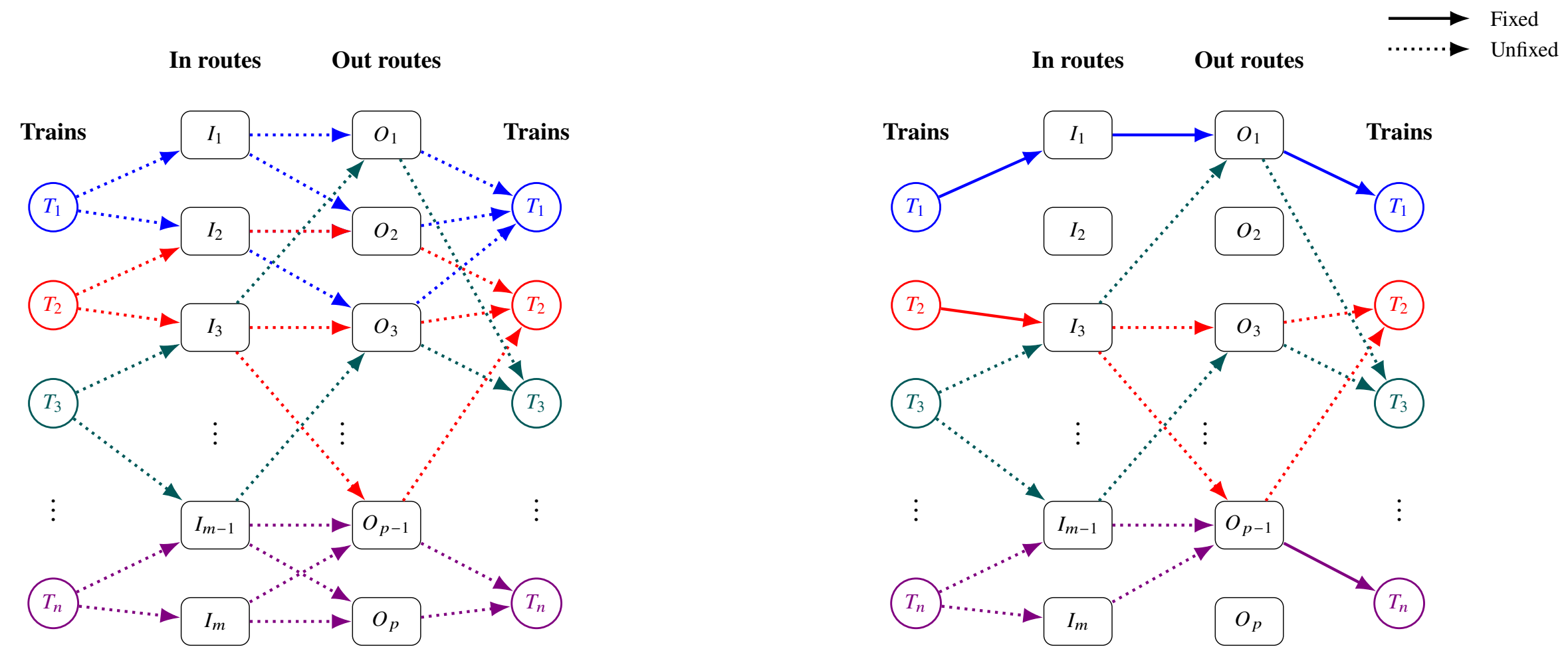
\begin{figure}[h]
    \centering
     \begin{subfigure}[t]{0.40\textwidth}
        \centering
        \resizebox{!}{6cm}{\begin{tikzpicture}[
    >=Latex,
    node distance=0.7cm and 2.0cm,
    train/.style={circle, draw, thick, minimum size=7mm, font=\small},
    route/.style={rectangle, draw, rounded corners, minimum width=10mm, minimum height=7mm, font=\small},
    dots/.style={font=\large},
    edgeT1/.style={->, very thick, dotted, blue},
    edgeT2/.style={->, very thick, dotted, red},
    edgeT3/.style={->, very thick, dotted, teal!70!black},
    edgeTn/.style={->, very thick, dotted, violet}
]

\node[train, draw=blue, text=blue] (T1s) {$T_1$};
\node[train, draw=red, text=red, below=of T1s] (T2s) {$T_2$};
\node[train, draw=teal!70!black, text=teal!70!black, below=of T2s] (T3s) {$T_3$};
\node[dots, below=of T3s] (TdotsS) {$\vdots$};
\node[train, draw=violet, text=violet, below=of TdotsS] (Tns) {$T_n$};

\node[route, right=1.5cm of T2s, yshift=2.5cm] (In1) {$I_{1}$};
\node[route, below=of In1] (In2) {$I_2$};
\node[route, below=of In2] (In3) {$I_3$};
\node[dots, below=of In3] (Indots) {$\vdots$};
\node[route, below=of Indots] (Inm-1) {$I_{m-1}$};
\node[route, below=of Inm-1] (Inm) {$I_m$};

\node[route, right=1.5cm of In1] (Out1) {$O_1$};
\node[route, right=1.5cm of In2] (Out2) {$O_2$};
\node[route, right=1.5cm of In3] (Out3) {$O_3$};
\node[dots, right=1.5cm of Indots] (Outdots) {$\vdots$};
\node[route, right=1.5cm of Inm-1] (Outp-1) {$O_{p-1}$};
\node[route, below=of Outp-1] (Outp) {$O_p$};

\coordinate (EndCol) at ($(Out1)+(2.2cm,0)$);

\node[train, draw=blue, text=blue] 
    (T1e) at (EndCol |- T1s) {$T_1$};

\node[train, draw=red, text=red] 
    (T2e) at (EndCol |- T2s) {$T_2$};

\node[train, draw=teal!70!black, text=teal!70!black] 
    (T3e) at (EndCol |- T3s) {$T_3$};

\node[dots] 
    (TdotsE) at (EndCol |- TdotsS) {$\vdots$};

\node[train, draw=violet, text=violet] 
    (Tne) at (EndCol |- Tns) {$T_n$};

\draw[edgeT1] (T1s) -- (In1);
\draw[edgeT1] (T1s) -- (In2);

\draw[edgeT1] (In1) -- (Out1);
\draw[edgeT1] (In1) -- (Out2);
\draw[edgeT1] (In2) -- (Out2);
\draw[edgeT1] (In2) -- (Out3);

\draw[edgeT1] (Out1) -- (T1e);
\draw[edgeT1] (Out2) -- (T1e);
\draw[edgeT1] (Out3) -- (T1e);

\draw[edgeT2] (T2s) -- (In2);
\draw[edgeT2] (T2s) -- (In3);

\draw[edgeT2] (In2) -- (Out2);
\draw[edgeT2] (In3) -- (Out3);
\draw[edgeT2] (In3) -- (Outp-1);

\draw[edgeT2] (Out2) -- (T2e);
\draw[edgeT2] (Out3) -- (T2e);
\draw[edgeT2] (Outp-1) -- (T2e);

\draw[edgeT3] (T3s) -- (In3);
\draw[edgeT3] (T3s) -- (Inm-1);

\draw[edgeT3] (In3) -- (Out1);
\draw[edgeT3] (Inm-1) -- (Out3);

\draw[edgeT3] (Out1) -- (T3e);
\draw[edgeT3] (Out3) -- (T3e);

\draw[edgeTn] (Tns) -- (Inm-1);
\draw[edgeTn] (Tns) -- (Inm);
\draw[edgeTn] (Inm-1) -- (Outp-1);
\draw[edgeTn] (Inm-1) -- (Outp);
\draw[edgeTn] (Inm) -- (Outp-1);
\draw[edgeTn] (Inm) -- (Outp);
\draw[edgeTn] (Outp-1) -- (Tne);
\draw[edgeTn] (Outp) -- (Tne);

\node[font=\bfseries, above=0.5cm of T1s] {Trains};
\node[font=\bfseries, above=0.5cm of In1] {In routes};
\node[font=\bfseries, above=0.5cm of Out1] {Out routes};
\node[font=\bfseries, above=0.5cm of T1e] {Trains};

\end{tikzpicture}}
        \caption{Original}
     \end{subfigure}
    \hfill
    \begin{subfigure}[t]{0.56\textwidth}
        \centering
        \resizebox{!}{6cm}{\begin{tikzpicture}[
    >=Latex,
    node distance=0.7cm and 2.0cm,
    train/.style={circle, draw, thick, minimum size=7mm, font=\small},
    route/.style={rectangle, draw, rounded corners, minimum width=10mm, minimum height=7mm, font=\small},
    dots/.style={font=\large},
    edgeT1/.style={->, very thick, blue},
    edgeT2solid/.style={->, very thick, red},
    edgeT2dotted/.style={->, very thick, dotted, red},
    edgeT3solid/.style={->, very thick, teal!70!black},
    edgeT3dotted/.style={->, very thick, dotted, teal!70!black},
    edgeTnsolid/.style={->, very thick, violet},
    edgeTndotted/.style={->, very thick, dotted, violet}
]

\node[train, draw=blue, text=blue] (T1s) {$T_1$};
\node[train, draw=red, text=red, below=of T1s] (T2s) {$T_2$};
\node[train, draw=teal!70!black, text=teal!70!black, below=of T2s] (T3s) {$T_3$};
\node[dots, below=of T3s] (TdotsS) {$\vdots$};
\node[train, draw=violet, text=violet, below=of TdotsS] (Tns) {$T_n$};

\node[route, right=1.5cm of T2s, yshift=2.5cm] (In1) {$I_{1}$};
\node[route, below=of In1] (In2) {$I_2$};
\node[route, below=of In2] (In3) {$I_3$};
\node[dots, below=of In3] (Indots) {$\vdots$};
\node[route, below=of Indots] (Inm-1) {$I_{m-1}$};
\node[route, below=of Inm-1] (Inm) {$I_m$};

\node[route, right=1.5cm of In1] (Out1) {$O_1$};
\node[route, right=1.5cm of In2] (Out2) {$O_2$};
\node[route, right=1.5cm of In3] (Out3) {$O_3$};
\node[dots, right=1.5cm of Indots] (Outdots) {$\vdots$};
\node[route, right=1.5cm of Inm-1] (Outp-1) {$O_{p-1}$};
\node[route, below=of Outp-1] (Outp) {$O_p$};

\coordinate (EndCol) at ($(Out1)+(2.2cm,0)$);

\node[train, draw=blue, text=blue] 
    (T1e) at (EndCol |- T1s) {$T_1$};

\node[train, draw=red, text=red] 
    (T2e) at (EndCol |- T2s) {$T_2$};

\node[train, draw=teal!70!black, text=teal!70!black] 
    (T3e) at (EndCol |- T3s) {$T_3$};

\node[dots] 
    (TdotsE) at (EndCol |- TdotsS) {$\vdots$};

\node[train, draw=violet, text=violet] 
    (Tne) at (EndCol |- Tns) {$T_n$};

\draw[edgeT1] (T1s) -- (In1);

\draw[edgeT1] (In1) -- (Out1);

\draw[edgeT1] (Out1) -- (T1e);

\draw[edgeT2solid] (T2s) -- (In3);

\draw[edgeT2dotted] (In3) -- (Out3);
\draw[edgeT2dotted] (In3) -- (Outp-1);

\draw[edgeT2dotted] (Out3) -- (T2e);
\draw[edgeT2dotted] (Outp-1) -- (T2e);

\draw[edgeT3dotted] (T3s) -- (In3);
\draw[edgeT3dotted] (T3s) -- (Inm-1);

\draw[edgeT3dotted] (In3) -- (Out1);
\draw[edgeT3dotted] (Inm-1) -- (Out3);

\draw[edgeT3dotted] (Out1) -- (T3e);
\draw[edgeT3dotted] (Out3) -- (T3e);

\draw[edgeTndotted] (Tns) -- (Inm-1);
\draw[edgeTndotted] (Tns) -- (Inm);
\draw[edgeTndotted] (Inm-1) -- (Outp-1);
\draw[edgeTndotted] (Inm) -- (Outp-1);
\draw[edgeTnsolid] (Outp-1) -- (Tne);
\node[
    overlay,
    anchor=north east,
    fill=white,
    inner sep=4pt,
    font=\small
] at ($(current bounding box.north east)+(2.25cm,1.75cm)$) {
    \begin{tikzpicture}[baseline=(current bounding box.center)]
        \draw[->, very thick, black] (-0.4,0) -- (0.8,0);
        \node[anchor=west] at (0.95,0) {Fixed};

        \draw[->, very thick, dotted, black] (-0.4,-0.45) -- (0.8,-0.45);
        \node[anchor=west] at (0.95,-0.45) {Unfixed};
    \end{tikzpicture}
};

\node[font=\bfseries, above=0.5cm of T1s] {Trains};
\node[font=\bfseries, above=0.5cm of In1] {In routes};
\node[font=\bfseries, above=0.5cm of Out1] {Out routes};
\node[font=\bfseries, above=0.5cm of T1e] {Trains};

\end{tikzpicture}}
        \caption{Reduced}
     \end{subfigure}
    \caption{Impact of full ($T_1$) and partial ($T_2, T_n$) ML-induced route fixing on the routing search space.}
    \label{fig:route_graph}
\end{figure}

Routes are classified as stable if $\hat{y}_{j'r}=1$ and as potentially affected otherwise. The route-level predictions are then mapped to train-level routing variables ($\alpha,\beta$) through the nominal timetable. If the nominal route of a train is predicted as stable, the corresponding inbound or outbound route-selection binaries are fixed to their nominal values. Variables associated with routes predicted as potentially affected remain unfixed and optimizable. Consequently, the routing search space is reduced by partially or fully restricting the routing alternatives available to individual trains. \Cref{fig:route_graph} illustrates this reduction through a train--route compatibility graph, in which edges represent admissible train--inbound-route assignments and compatible inbound--outbound route pairs, while directed paths represent complete routing alternatives. To formalize these fixings, let $\mathcal{F}^{\mathrm{ML}}\subseteq B^r$ denote the set of route-selection binaries fixed by the Scoper, and let $\hat{u}_i\in\{0,1\}$ denote the corresponding nominal assignment for each $i\in\mathcal{F}^{\mathrm{ML}}$. The reduced rTPP is obtained by adding these fixing constraints to the original formulation:
\begingroup
\setlength{\abovedisplayskip}{4pt}
\setlength{\belowdisplayskip}{4pt}
\setlength{\belowdisplayshortskip}{4pt}
\begin{align}
\text{($\mathrm{rTPP}_{\mathrm{red}}$)}: \quad &\eqref{eq:trp_obj}\text{--}\eqref{eq:trp_cont}\nonumber\\
&u_i = \hat{u}_i \quad \forall i \in \mathcal{F}^{\mathrm{ML}}\label{eq7}
\end{align}
\endgroup
Since $\mathcal{F}^{\mathrm{ML}}$ contains only route-selection binaries, $\mathrm{rTPP}_{\mathrm{red}}$ reduces the routing search space while preserving the original objective, operational constraints, and domains of all unfixed variables. The unfixed routing variables and the remaining variables $v=(d, q)$ are therefore jointly optimized in the resulting restricted MILP. Because $\mathcal{F}^{\mathrm{ML}}$ is derived from route-stability predictions, two types of prediction error may arise. A false negative (FN) classifies a stable route as potentially affected and leaves the associated routing variable unfixed, thereby forgoing part of the attainable search-space restriction. A false positive (FP), in contrast, classifies an affected route as stable and imposes route fixings that may exclude routing alternatives needed in the rescheduled solution. Because waiting time, dwell time, and delay extensions are bounded in realistic railway operations, an incorrect fixing may induce route, platform, or headway conflicts that cannot be resolved by timing adjustments alone. Consequently, $\mathrm{rTPP}_{\mathrm{red}}$ may be infeasible even when the original rTPP remains feasible.

The predictions could alternatively be imposed as soft restrictions through an additional objective penalty. We retain hard route fixings because the original rTPP already penalizes deviation from the nominal routes through $\lambda_{\mathcal{R}}$ (see \hyperref[app:rtpp_formulation]{Appendix}), while the learned fixings are intended to restrict the corresponding routing decisions. An additional soft penalty would therefore penalize the same route deviation twice, causing its effective cost to depend on the predicted route status and requiring calibration against the existing cost components. Hard route fixings provide an explicit search-space reduction without altering the original objective, with feasibility restored by repair and solution quality addressed by refinement.

\subsection{Feasibility repair}\label{sec:repair}
The repair stage restores feasibility by selectively relaxing the ML-induced route-fixing constraints indexed by $\mathcal{F}^{\mathrm{ML}}$ in $\mathrm{rTPP}_{\mathrm{red}}$, while preserving all original operational constraints of the rTPP. We consider two repair strategies. The first is an iterative conflict-guided repair, which uses solver-identified conflicts to restrict the candidate fixings for relaxation. The second is a globally scoped repair, which considers all ML-induced route fixings simultaneously and is implemented through an increasing sequence of cardinality-restricted MILPs.

\subsubsection{Conflict-guided repair}
The conflict-guided repair restricts each repair MILP to fixings implicated in the current infeasibility. Let $\mathcal{F}^{k}\subseteq\mathcal{F}^{\mathrm{ML}}$ denote the set of active ML-induced route fixings at repair iteration $k$, with $\mathcal{F}^{0}=\mathcal{F}^{\mathrm{ML}}$. The corresponding reduced model $\mathrm{rTPP}_{\mathrm{red}}^{k}$ enforces $u_i=\hat{u}_i$ for all $i\in\mathcal{F}^{k}$. When $\mathrm{rTPP}_{\mathrm{red}}^{k}$ is infeasible, the solver's conflict refiner identifies a minimal set of mutually incompatible constraints and variable bounds. Let $\mathcal{J}^{k}$ denote the original rTPP constraints in this conflict, each represented as $g_j(u,v),\bowtie_j,b_j$, where $\bowtie_j\in\{\leq,=,\geq\}$. 

Let $\mathcal{Q}^{k}\subseteq \mathcal{F}^{k}$ denote the active fixings associated with routing variables that appear in at least one constraint of $\mathcal{J}^{k}$. Because the constraints in $\mathcal{J}^{k}$ and the fixings in $\mathcal{Q}^{k}$ are jointly infeasible, releasing an arbitrary fixing from $\mathcal{Q}^{k}$ need not resolve the conflict. The auxiliary repair MILP therefore 
considers only $\mathcal{Q}^{k}$ as release candidates and replaces their hard fixing constraints with release-controlled constraints controlled by binary release variables $\rho_i$:
\begingroup
\setlength{\abovedisplayskip}{4pt}
\setlength{\belowdisplayskip}{4pt}
\setlength{\belowdisplayshortskip}{4pt}
\begin{align}
(\mathrm{Repair}_{\mathrm{CG}}^{k}) \quad 
\min_{u,v,\rho} \quad & \sum_{i\in \mathcal{Q}^k}\rho_i \label{eq:repair_obj1}\\
\text{s.t.} \quad 
& g_j(u,v)\bowtie_j,b_j & \forall j\in \mathcal{J}^k\label{eq:iis_sense}\\
& u_i \leq \hat{u}_i + \rho_i &\forall i\in \mathcal{Q}^k\label{eq:iis_link1}\\
& u_i \geq \hat{u}_i - \rho_i &\forall i\in \mathcal{Q}^k\label{eq:iis_link2}\\
&\rho_i \in \{0, 1\} &\forall i\in \mathcal{Q}^k\label{eq:iis_binary}\\
& (u,v)\in\mathcal{D} &\label{eq:iis_domains}
\end{align}
\endgroup
At repair iteration $k$, the objective \eqref{eq:repair_obj1} minimizes the number of released fixings in $\mathcal{Q}^{k}$, thereby identifying a minimum-cardinality release that renders the current conflict subsystem satisfiable. For each $i\in\mathcal{Q}^{k}$, the binary variable $\rho_i$ determines whether the prescribed value of the route-selection variable $u_i$ is retained or released. Constraints \eqref{eq:iis_sense} retain the original senses of the rTPP constraints contained in the conflict, while $\mathcal{D}$ denotes the variable domains defined in \eqref{eq:trp_routebin}\text{--}\eqref{eq:trp_cont}. Constraints \eqref{eq:iis_link1} and \eqref{eq:iis_link2} replace the hard fixing $u_i=\hat{u}_i$ with a release-controlled relation. If $\rho_i=0$, the fixing is enforced; if $\rho_i=1$, it is released and $u_i$ may take either value in its original binary domain. These constraints are used only in the auxiliary repair MILP. Because its objective contains no operational cost terms, the release variables neither alter the rTPP objective nor require calibration against its existing cost components. Only constraints associated with ML-induced fixings are eligible for release; no original rTPP constraint is relaxed.

After solving ($\mathrm{Repair}_{\mathrm{CG}}^{k}$), let $\rho^{k,*}$ denote the optimal release vector. The active fixing set is updated as $\mathcal{F}^{k+1} = \mathcal{F}^{k}\setminus\{i\in\mathcal{Q}^{k}:\rho_i^{k,*}=1\}$. Thus, only fixings selected from $\mathcal{Q}^{k}$ are removed; all remaining fixings are enforced in the subsequent model $\mathrm{rTPP}{\mathrm{red}}^{k+1}$. The resulting reduced model $\mathrm{rTPP}_{\mathrm{red}}^{k+1}$ is then re-examined for feasibility. Conflict identification and repair are repeated until feasibility is restored or the prescribed iteration limit is reached. Each successful repair iteration releases at least one fixing, so the active fixing set strictly decreases. Provided that every remaining infeasibility yields a nonempty candidate set $\mathcal{Q}^{k}$, the procedure requires at most $|\mathcal{F}^{\mathrm{ML}}|$ repair iterations. Upon successful termination, the remaining active fixing set is denoted by $\mathcal{F}^{\mathrm{R}}$.

\subsubsection{Global repair model}
The conflict-guided repair restricts the release candidate to the current conflict and may require several iterations when infeasibility involves multiple interacting fixings. In contrast, the global repair considers all ML-induced route fixings in $\mathcal{F}^{\mathrm{ML}}$ simultaneously and identifies the minimum number that must be released to recover feasibility. This is measured by the Hamming distance between the routing vector $u$ and the prescribed values $\hat{u}$ over $\mathcal{F}^\mathrm{ML}$: 
\begin{equation}
D_{\mathrm H}(u,\hat{u})
=
\sum_{\substack{i\in\mathcal{F}^{\mathrm{ML}}\\ \hat{u}_i=0}} u_i
+
\sum_{\substack{i\in\mathcal{F}^{\mathrm{ML}}\\ \hat{u}_i=1}} (1-u_i).
\label{eq:hamming}
\end{equation}
Because $u_i, \hat{u}_i\in \{0,1\}$, each term equals one when $u_i$ differs from its prescribed value and zero otherwise. The global repair model is formulated as: 

\begin{equation}
\label{eq:repair_obj}
\begin{aligned}
(\mathrm{Repair}_{\mathrm{G}}) \quad 
\min_{u,v} \quad & D_{\mathrm H}(u,\hat{u}) \\
\text{s.t.} \quad & \eqref{eq:trp_time}\text{--}\eqref{eq:trp_cont}
\end{aligned}
\end{equation}

Model ($\mathrm{Repair}_{\mathrm{G}}$) minimizes the number of deviations from the prescribed values while enforcing all operational constraints and variable domains of the original rTPP. Since $u_i$ and $\hat{u}_i$ are binary, the Hamming-distance expression in \eqref{eq:hamming} is linear and requires no auxiliary release variables; it corresponds to the distance function used in the feasibility pump \citep{fischetti2005feasibility}. Because the rTPP contains big-$M$ coefficients in its disjunctive sequencing and headway constraints (see \hyperref[app:rtpp_formulation]{Appendix}), we solve $\mathrm{Repair}{\mathrm{G}}$ as an integral MILP rather than apply a feasibility-pump scheme, whose alternation between LP optimization and rounding can exhibit pathological behavior in such formulations \citep{fischetti2005feasibility}. For an optimal repair solution $(u^{\mathrm R},v^{\mathrm R})$, the retained fixing set is $\mathcal{F}^{\mathrm R} =\{i\in \mathcal{F}^{\mathrm{ML}}: u_i^{\mathrm R}=\hat{u}_i\}$, while fixings satisfying $u_i^{\mathrm R}\neq\hat{u}_i$ are released. Thus, all original rTPP constraints remain enforced, and only the ML-induced routing fixings may be released. Because feasibility is expected to require only a small number of releases, we solve the global repair through an increasing sequence of cardinality-restricted problems:
\begingroup
\setlength{\abovedisplayskip}{4pt}
\setlength{\belowdisplayskip}{4pt}
\setlength{\belowdisplayshortskip}{4pt}
\begin{align}
(\mathrm{Repair}_{\mathrm G}^{K}) \qquad
\min_{u,v}\quad
& D_{\mathrm H}(u,\hat{u})
\label{eq:repair_k_obj}\\
\mathrm{s.t.}\quad
& \eqref{eq:trp_time}\text{--}\eqref{eq:trp_cont},
\nonumber\\
& D_{\mathrm H}(u,\hat{u})\leq K.
\label{eq:cardinality}
\end{align}
\endgroup
The bound $K$ limits the number of fixings that may be released, while all fixings in $\mathcal{F}^{\mathrm{ML}}$ remain eligible for release. Constraint \eqref{eq:cardinality} has the form of a local branching constraint around $\hat{u}$, restricted here to the ML-induced fixings \citep{fischetti2003local}. Starting from a small value, $K$ is increased until $\mathrm{Repair}_{\mathrm{G}}^{K}$ admits a feasible solution, indicating that feasibility can be restored by releasing at most $K$ fixings. Solving the first feasible model in this sequence to optimality yields the same minimum Hamming distance as $\mathrm{Repair}{\mathrm{G}}$. The $K$-values therefore need not be consecutive; they only control the successive sizes of the restricted repair region.

The global repair is a feasibility-restoration model rather than a rescheduling model. Its solution determines the repaired fixing set $\mathcal{F}^{\mathrm{R}}$ and provides a feasible starting solution for the subsequent refinement and optimization stages, in which the original rTPP objective is reinstated.

\subsection{Interaction-guided refinement}\label{sec:refine}
The refinement stage recovers routing alternatives that may support high-quality feasible solutions but remain excluded after reduction and, when required, repair. Rather than releasing route fixings uniformly, refinement selectively relaxes the restricted MILP by releasing fixings associated with trains whose cancellations or positive waiting- or dwell-time extensions contribute to the incumbent rTPP objective value. The release set is extended to their interaction neighbors, since reducing these deviations may require coordinated routing changes among trains with shared infrastructure and operate in temporal proximity. Train interactions are represented by an interaction graph $\mathcal{G}^{\mathrm{int}}=(\mathcal{T},E^{\mathrm{int}})$, where each node represents a train and an edge $(t,\dot{t})\in E^{\mathrm{int}}$ indicates a potential interaction arising from spatial overlap between the trains' feasible routes in $\mathcal{G}^{\mathcal{I}}$ and temporal proximity in the nominal timetable. The graph thus provides a compact representation of train pairs whose routing decisions are likely to be coupled.

Let $(u^{\mathrm{R}},v^{\mathrm{R}})$ denote the feasible incumbent available when refinement begins, obtained directly from the reduced model or after feasibility repair. Let $\mathcal{F}^{\mathrm{R}}\subseteq\mathcal{F}^{\mathrm{ML}}$ denote the route fixings that remain active, with $\mathcal{F}^{\mathrm{R}}=\mathcal{F}^{\mathrm{ML}}$ when no repair is required. The mapping $\tau:\mathcal{F}^{\mathrm{R}}\rightarrow \mathcal{T}$ assigns each fixing $i\in\mathcal{F}^{\mathrm{R}}$ to its associated train. Let $\mathcal{T}^{\mathrm{aff}}$ denote the trains that are canceled or have a positive waiting- or dwell-time extension relative to the nominal timetable in $(u^{\mathrm{R}},v^{\mathrm{R}})$. Their interaction neighborhood is $N(\mathcal{T}^{\mathrm{aff}})=\{\dot{t}\in \mathcal{T}:(t,\dot{t})\in E^{\mathrm{int}} \text{ for some } t\in \mathcal{T}^{\mathrm{aff}}\}$. The set of fixings selected for release is, $\mathcal{H}=\{i\in\mathcal{F}^{\mathrm{R}}:
\tau(i)\in\mathcal{T}^{\mathrm{aff}}\cup N(\mathcal{T}^{\mathrm{aff}})\}$. The refined model retains all remaining fixings:
\begingroup
\setlength{\abovedisplayskip}{4pt}
\setlength{\belowdisplayskip}{4pt}
\setlength{\belowdisplayshortskip}{4pt}
\begin{align}
(\mathrm{rTPP}_{\mathrm{ref}}):\quad 
& \eqref{eq:trp_obj}\text{--}\eqref{eq:trp_cont}\nonumber\\
&u_i=\hat{u}_i  &\forall i\in\mathcal{F}^{\mathrm{R}}\setminus\mathcal{H}.
\end{align}
\endgroup
Routing binaries indexed by $\mathcal{H}$ recover their original binary domains and are re-optimized, while all other active fixings remain enforced. Consequently, $\mathrm{rTPP}_{\mathrm{ref}}$ is a relaxation of the restricted model defined by  $\mathcal{F}^{\mathrm{R}}$ and remains a restriction of the original rTPP. Because refinement removes only route-fixing constraints, the incumbent $(u^{\mathrm{R}},v^{\mathrm{R}})$ remains feasible for $\mathrm{rTPP}_{\mathrm{ref}}$, and its objective value provides an upper bound on the optimal value of the refined model.

\subsection{Optimization}\label{sec:optimization}
The optimization stage solves the current restricted rTPP under the original rTPP objective. The solver jointly optimizes all unfixed routing binaries and the remaining decision variables, while enforcing the retained route fixings. These fixings restrict the admissible routing assignments and concentrate the branch-and-bound search on the region retained by the preceding stages. Because the restricted model differs from the original rTPP only through additional route-fixing constraints, every feasible solution of the restricted model is also feasible for the original rTPP.

\subsection{Online execution of 3RO}\label{sec:online_3RO}
Let $\mathrm{rTPP}(\mathcal{F})$ denote the original rTPP augmented with the fixing constraints $u_i=\hat{u}_i$ for all $i\in\mathcal{F}$, where $\mathcal{F}$ is the current set of active route fixings. \Cref{alg:online_3ro} summarizes the online execution of 3RO under a common wall-clock budget $B$. Applying the trained Scoper to the disrupted instance produces the ML-induced route fixings $\mathcal{F}^{\mathrm{ML}}$. If $\mathrm{rTPP}(\mathcal{F}^{\mathrm{ML}})$ is feasible, we set $\mathcal{F}^{\mathrm{R}}=\mathcal{F}^{\mathrm{ML}}$; otherwise, the selected repair strategy seeks a retained set $\mathcal{F}^{\mathrm{R}}\subseteq\mathcal{F}^{\mathrm{ML}}$ for which $\mathrm{rTPP}(\mathcal{F}^{\mathrm{R}})$ is feasible. Let $(u^{\mathrm{inc}},v^{\mathrm{inc}})$ denote the best feasible incumbent maintained during online execution. Once feasibility of $\mathrm{rTPP}(\mathcal{F}^{\mathrm{R}})$ has been established, the model is solved under the original rTPP objective to obtain the initial incumbent used by refinement.

Starting from $\mathcal{F}\leftarrow\mathcal{F}^{\mathrm{R}}$, refinement alternates selective fixing release and optimization. At each iteration, the interaction-neighborhood rule uses the current incumbent to determine a release set $\mathcal{H}\subseteq\mathcal{F}$. If $\mathcal{H}\neq\emptyset$, the active set is updated as $\mathcal{F}\leftarrow\mathcal{F}\setminus\mathcal{H}$ and the resulting restricted model is optimized; the updated incumbent then guides the next refinement iteration. If $\mathcal{H}=\emptyset$, no additional fixings are released. Consequently, each refined model is a relaxation of its predecessor, and the corresponding feasible regions expand monotonically. Because each continuing iteration releases at least one active fixing, at most $|\mathcal{F}^{\mathrm{R}}|$ fixing-release iterations are possible. This common active-set representation links repair, which determines a retained fixing set, with refinement, which selectively releases additional fixings to recover routing alternatives.

The nested relaxations also permit incumbent transfer: removing fixing constraints preserves the feasibility of the current incumbent, which can therefore be supplied as a warm start to the subsequent model. For this minimization problem, retaining the best incumbent across successive optimization calls yields a nonincreasing sequence of incumbent objective values, although an individual time-limited solve need not produce an improvement. 

Let $\ell$ denote the elapsed wall-clock time from the start of online execution, and let $b_{\mathrm{init}}$ and $b_{\mathrm{ref}}$ denote the solve-time limits for the initial restricted model and each refined model, respectively. All online stages contribute to $\ell$, and every solver call is capped by the remaining time $B-\ell$, ensuring that the total online computation does not exceed $B$. The initial restricted model is solved for at most $\min\{b_{\mathrm{init}},B-\ell\}$, and each refined model for at most $\min\{b_{\mathrm{ref}},B-\ell\}$. If $\mathcal{H}=\emptyset$, the current restricted model is solved for up to the remaining time $B-\ell$.

If repair fails to restore feasibility, all active route fixings are released and the unrestricted rTPP is solved using the remaining budget. If no feasible incumbent is obtained within the initial solve-time limit, the same restricted model continues to be solved using the remaining time. Once a feasible incumbent is available, the procedure operates as an anytime method and returns the best incumbent available at termination.

\begin{algorithm}[t]
\caption{Online execution of 3RO under a fixed wall-clock budget}
\label{alg:online_3ro}
\begin{algorithmic}[1]
\REQUIRE Disrupted instance $j'$; trained Scoper; time-budget parameters: $B$, $b_{\mathrm{init}}$, and $b_{\mathrm{ref}}$
\ENSURE Best feasible incumbent $(u^{\mathrm{inc}},v^{\mathrm{inc}})$ found within $B$

\STATE Start the timer; let $\ell$ denote the elapsed time and initialize with no feasible incumbent
\STATE Apply the trained Scoper to $j'$, construct $\mathcal{F}^{\mathrm{ML}}$ and set $\mathcal{F}\leftarrow\mathcal{F}^{\mathrm{ML}}$
\COMMENT{Reduction; \Cref{sec:reduction}}

\IF{$\mathrm{rTPP}(\mathcal{F})$ is infeasible}
    \STATE Apply the selected repair strategy within $B-\ell$
    \COMMENT {Repair; \Cref{sec:repair}}
    \IF{repair fails to restore feasibility}
    \STATE Set $\mathcal{F}\leftarrow\emptyset$ and solve $\mathrm{rTPP}(\mathcal{F})$ for $B-\ell$; update and \textbf{return} the incumbent
    \ENDIF
    \STATE Set $\mathcal{F}\leftarrow\mathcal{F}^{\mathrm{R}}$
\ENDIF
    
\STATE Solve $\mathrm{rTPP}(\mathcal{F})$ for $\min\{b_{\mathrm{init}},B-\ell\}$; update the incumbent
\COMMENT{Optimize; \Cref{sec:optimization}}
\IF{no feasible incumbent is available}
\STATE Continue solving $\mathrm{rTPP}(\mathcal{F})$ for $B-\ell$; update \textbf{return} the incumbent
\ENDIF

\WHILE{$\ell<B$}
\STATE Determine $\mathcal{H}\subseteq\mathcal{F}$ from the current incumbent
\COMMENT{Refinement; \Cref{sec:refine}}
\IF{$\mathcal{H}=\emptyset$}
\STATE Continue solving $\mathrm{rTPP}(\mathcal{F})$ for $B-\ell$; update and \textbf{return} the incumbent
\ENDIF
\STATE Set $\mathcal{F}\leftarrow\mathcal{F}\setminus\mathcal{H}$
\STATE Solve $\mathrm{rTPP}(\mathcal{F})$, warm-started with the incumbent, for $\min{b_{\mathrm{ref}},B-\ell}$; update the incumbent
\ENDWHILE

\STATE \textbf{return} $(u^{\mathrm{inc}},v^{\mathrm{inc}})$
\end{algorithmic}
\end{algorithm}

\section{Computational setup}
\label{sec:computational_study}
This section describes the computational setup used to evaluate the proposed 3RO framework. We first introduce the case-study setting, including station infrastructure, nominal timetable, and disruption scenarios, followed by generation of the learning dataset and its partitioning into training, validation, and evaluation sets. The evaluation set is further stratified according to the computational effort required by the unrestricted MILP, distinguishing instances that are tractable within the real-time horizon from those requiring substantially longer solution times. The experimental protocol specifies the comparison of 3RO with the full CPLEX MILP under common real-time and extended computation limits. Finally, we summarize the implementation environment and solver settings used throughout the experiments.

\subsection{Case-study and disruption scenarios}
\label{sec:case_study_setting}
\paragraph{Infrastructure and timetable.}
The computational study considers a medium-sized railway station in the German rail network with eight platform tracks and $|\mathcal{R}|= 125$ predefined routes, of which $|\mathcal{R}^{\mathrm{nom}}|=72$  are used in the nominal timetable and remaining $53$ serve as rerouting alternatives. The timetable spans a six-hour operating period (05:00--11:00) with more than $250$ train movements and peak traffic of approximately $45$ trains per hour. Each train $t\in\mathcal{T}$ is characterized by its scheduled arrival and departure times, planned platform track, nominal inbound--outbound route pair, and its feasible route alternative $\mathcal{R}^{+}_t$ and $\mathcal{R}^{-}_t$, which contain $12$--$15$ alternatives per direction on average.

\paragraph{Disruption scenarios.}
We consider two classes of capacity-reducing disruptions in the station area: switch failures and platform closures. A switch failure renders all routes traversing the affected switch unavailable for the duration of the disruption, whereas a platform closure removes a platform track together with its associated inbound and outbound routes. Both disruption types reduce the available routing and platforming capacity and intensify interactions between trains arising from route incompatibility, platform occupancy, and minimum headway requirements. However, as they differ in how the resulting capacity loss is distributed over the station topology, they are therefore retained as separate classes.

A disruption instance is characterized by its class, location, onset time, and duration. Candidate locations cover approximately 20\% of the station infrastructure, with switch failures sampled three times as frequently as platform closures. Disruption durations range from $15$ to $120$ minutes in $5$-minute increments, following a right-skewed distribution that favors longer disruptions, which impose higher computational burden \citep{zhan2015real}. Onset times are sampled across the operating period with greater representation of peak traffic. The resulting scenario set therefore captures variation in disruption location, timing, and severity.

\subsection{Learning and evaluation data}
\paragraph{Reference solutions and learning data.}
Each generated disruption instance is solved with the original rTPP MILP under a $4\%$ relative optimality-gap tolerance and a one-hour time limit, using fixed solver settings and a single thread to ensure reproducible reference solutions\footnote{these thresholds are chosen only to keep the computation times manageable.}. Each instance-solution pair forms one learning sample. The input features follow the construction described in \Cref{sec:reduction} and are derived exclusively from the disruption declaration and nominal timetable. The target vector contains one route-stability label for each of the $|\mathcal{R}^{\mathrm{nom}}|=72$ nominal routes, yielding a fixed-dimensional multi-label target. Since only reference solutions obtained within the one-hour time limit and satisfying the $4\%$ relative optimality-gap tolerance are retained, the route-stability labels are derived from near-optimal routing plans that are not necessarily proven optimal. Moreover, as discussed in \Cref{sec:reduction}, route-stability labels are defined relative to the selected reference solution and may therefore differ when alternative optimal routing plans exist.

The dataset comprises $2,412$ instances, partitioned into $1,586$ training, $600$ validation, and $226$ evaluation instances. The evaluation set is excluded from all model fitting and selection and contains previously unobserved combinations of disruption location, onset time, and duration. This provides an out-of-sample assessment of generalization to unseen disruption scenarios under spatial and temporal variation within the studied station.

\paragraph{Difficulty stratification.}
The $226$ evaluation instances are stratified by the computational effort the original rTPP MILP requires to reach a $4\%$ relative optimality gap: at most $180$~s (easy), more than $180$~s but at most $360$~s (medium), and more than $360$~s or not reaching the target within the one-hour run (hard). This yields $74$ easy, $24$ medium, and $128$ hard instances. The stratification distinguishes instances solved within the real-time horizon from those requiring progressively greater computational effort, enabling the performance of 3RO to be examined as instance difficulty increases.

\subsection{Experimental design}
\label{sec:experimental_design}
\paragraph{Solution approaches and time limits.} The computational study evaluates the proposed 3RO framework against three complementary reference approaches: the unrestricted MILP, a variant without refinement (2RO), and the greedy rTPP heuristic. The \emph{unrestricted MILP} solves the original rTPP without ML-induced route fixings and serves as the full-model benchmark. The proposed \emph{3RO} framework follows the reduce--repair--refine--optimize procedure described in \Cref{sec:methodology}. The \emph{2RO} variant omits interaction-guided refinement and directly optimizes the restricted model obtained after reduction and, when required, repair. The initially reduced model represents the learned restriction induced by the Scoper; accordingly, 2RO follows reduce-optimize when this model is feasible and reduce-repair-optimize otherwise. Comparing 3RO with 2RO therefore isolates the specific contribution of refinement. 

As an external rTPP benchmark, we implement the greedy heuristic of \citet{garcia2024efficient} and adapt it to the infrastructure representation and operational requirements of the present problem. The underlying sequential solution procedure is retained, while the implementation is extended where necessary to ensure compatibility with the feasible routing alternatives, operational constraints, and objective function used in our experiments. This allows the heuristic and 3RO solutions to be evaluated under the same feasibility and solution-quality criteria.

The MILP-based approaches are evaluated under a common wall-clock budget of $B=180$~s \citep{sama2016ant, pascariu2024formulation}. This budget represents the operational computation time available for obtaining a revised plan after the disruption is observed and covers the complete online phase defined in \Cref{sec:methodology}. For 2RO, the restricted model is optimized using the budget remaining after reduction and repair. Within 3RO, the initial restricted solve and each refined solve are limited to $b_{\mathrm{init}}=b_{\mathrm{ref}}=60$~s, subject to the remaining global budget. Extended full-MILP runs of up to 3,600~s provide a long-horizon solution-quality reference.

\paragraph{Scoper specifications.}
The candidate model set is $\mathcal{M}={\mathrm{LR},\mathrm{RF},\mathrm{MLP}}$, comprising logistic regression, random forest, and multilayer perceptron classifiers. For each route--model pair, the candidate thresholds are the default value $\phi^{\mathrm{def}}=0.5$ and the ROC-based threshold $\phi^{\mathrm{ROC}}_{r,m}$. Residual ties in the route-specific model--threshold selection are resolved using the fixed priority RF$\rightarrow$MLP$\rightarrow$LR. For comparison, homogeneous Scoper baselines apply a single model class across all route labels, using either the default threshold or ROC-derived route-specific thresholds.

\paragraph{Repair strategies and settings.}
For the global repair, the cardinality-restricted formulation is solved sequentially with $K=2,5$ and thereafter $K\leftarrow\min\{2K,|\mathcal{F}^{\mathrm{ML}}|\}$ until feasibility is restored. Because the Hamming-distance objective is integer-valued, an absolute MIP-gap tolerance strictly below one certifies optimality: once the incumbent--bound gap is below one, no better integer objective value can exist. The repair comparison includes the proposed conflict-guided and global strategies together with CPLEX FeasOpt as a solver-native benchmark. FeasOpt restores feasibility by minimizing a penalty associated with permitted relaxations of variable bounds or constraints. We allow only the bounds imposing ML-induced route fixings to be relaxed, assign them equal penalties, and retain all original rTPP constraints. Under minimum-sum mode, changing a fixed binary value requires a unit bound relaxation; hence, the FeasOpt objective corresponds to the number of released fixings.

\paragraph{Evaluation criteria.}
The analysis addresses five aspects of the framework. First, the \emph{learning performance} is evaluated through the route-level model--threshold pairs selected during validation and their out-of-sample macro-$F_1$ and macro-AP values relative to the homogeneous Scoper baselines. Macro-$F_1$ and macro-AP denote the unweighted mean across route labels of the per-route $F_1$ and average precision values, respectively, as introduced in \Cref{sec:reduction}. Second, the \emph{learning-based reduction} is assessed through the fixing rate $|\mathcal{F}^{\mathrm{ML}}|/|B^{r}|$, false-positive fixings, and the frequency of feasibility repair. Together, these measures characterize the trade-off between reducing the routing search space and preserving feasible routing alternatives. Third, the \emph{repair strategies} are compared in terms of repair time and the number of released fixings. 

Fourth, the contribution of \emph{refinement} is assessed by comparing 3RO with 2RO in terms of terminal solution quality and the timing of primal improvements. This comparison is further stratified according to whether the initially reduced model is feasible or requires repair, allowing the effect of refinement to be examined both after direct learned scope restriction and after feasibility restoration. Finally, the overall performance of 3RO is evaluated against the greedy rTPP heuristic and the unrestricted MILP. The comparison with the heuristic examines the evolution of feasibility and solution quality over the 180~s horizon, while the unrestricted-MILP comparison considers both the common 180~s horizon and the extended 3,600~s reference. Relative solution improvement is computed as $100\times(z^{\mathrm{ref}}-z)/max\{|z^{\mathrm{ref}}|,\epsilon\}$, where $z$ is the objective value of the evaluated solution approach, $z^{\mathrm{ref}}$ is the corresponding benchmark value, and $\epsilon>0$ prevents division by zero. Positive values indicate improvement over the reference.

\paragraph{Implementation}
The framework is implemented in Python~3.9, and all MILP models are solved using IBM ILOG CPLEX~22.1.1. Logistic regression and random forest models are implemented using scikit-learn \citep{pedregosa2011scikit}, while the multilayer perceptron is implemented in TensorFlow/Keras \citep{abadi2016tensorflow,chollet2015keras}. All experiments are conducted on a PC equipped with an AMD EPYC~7513 processor at 2.60\,GHz and 256\,GB of memory.

\section{Results}
This section presents the computational results for the proposed 3RO framework on a realistic rail station area. It evaluates how the learning-based reduction and the subsequent repair and refinement stages affect the computational performance and solution quality of 3RO. The analysis proceeds from the predictive behavior of the Scoper (\Cref{sec:results_scoper}) to its effects on the restricted MILP (\Cref{sec:results_downstream_learn,sec:results_downstream_repair,sec:results_downstream_refine}) and, finally, to the performance of the complete framework under real-time budget (\Cref{sec:results_sol_quality}).

\subsection{Route-level heterogeneity in prediction.}\label{sec:results_scoper}
This section examines route-level heterogeneity in the Scoper from two complementary perspectives. We first characterize the model--threshold pairs selected across route labels and then assess the predictive gains of route-specific selection relative to the homogeneous Scoper baselines defined in \Cref{sec:experimental_design}.

\noindent\textbf{Route-level configuration heterogeneity.}
We first examine whether the offline validation results support selecting the model class and threshold separately for each route label. \Cref{tab:het_composition} reports the composition of the selected model--threshold pairs across the $72$ nominal routes. All three model classes are represented, although RF is dominant, being selected for 63.9\% of the labels, compared with 22.2\% for LR and 13.9\% for MLP. Threshold selection is more concentrated, with the default threshold $\phi^{\mathrm{def}}=0.5$ selected for 80.6\% of the labels and the ROC-based threshold for the remaining 19.4\%. The most frequent pair, RF with $\phi^{\mathrm{def}}=0.5$, accounts for only 48.6\% of the labels. Thus, no single model--threshold pairs is selected for a majority of route labels, supporting the route-specific construction of the Scoper.

\begin{table}[h]
\centering
\caption{Composition of the selected model--threshold pairs across the
$|\mathcal{R}^{\mathrm{nom}}|=72$ route labels}
\label{tab:het_composition}
\renewcommand{\arraystretch}{1.5}
\begin{tabular}{lccc}
\hline
 & \multicolumn{2}{c}{Threshold} & \\
\cline{2-3}
Model class & $\phi^{\mathrm{def}}=0.5$ & $\phi^{\mathrm{ROC}}_{r,m}$ & Labels \\
\hline
RF  & 48.6 & 15.3 & 46 \\
LR  & 19.4 &  2.8 & 16 \\
MLP & 12.5 &  1.4 & 10 \\
\hline
Labels & 58 & 14 & 72 \\
\hline
\end{tabular}
\end{table}

\noindent\textbf{Predictive value of heterogeneous selection.}
We compare the heterogeneous Scoper with six homogeneous Scoper variants, obtained by applying one model class and one threshold rule across all route labels. \Cref{tab:het_vs_hom} reports the difference $\Delta=\mathrm{HET}-\mathrm{HOM}$; where HET and HOM denote the metric values obtained by the heterogeneous Scoper and the corresponding homogeneous variant, respectively. Positive values therefore favor the heterogeneous Scoper. The differences are positive for both Macro-$F_1$ and macro-AP against every homogeneous variant and within every evaluation subset, indicating that the advantage of route-specific selection persists across difficulty levels. 

RF provides the closest homogeneous alternative, with the smallest overall differences of $0.0047$ in Macro-$F_1$ and $0.0038$ in macro-AP, consistent with RF being selected for most route labels. The differences are substantially larger when LR or MLP is applied across all labels, reaching $0.1167$ in Macro-$F_1$ relative to LR with ROC-derived thresholds. Importantly, because macro-AP is threshold-independent, its consistently positive differences indicate that the advantage does not arise solely from route-specific threshold selection; selecting the model class by route label also improves probability ranking. Although the magnitude of the gain varies across homogeneous variants and difficulty levels, its direction remains consistent. Because the Scoper predictions determine the set of fixings $\mathcal{F}^{\mathrm{ML}}$, even small predictive differences can change which routing variables are fixed and, consequently, the feasible region of the reduced rTPP. The following analysis examines these downstream effects.

\begin{table}[h]
\centering
\caption{Difference in prediction performance between the heterogeneous Scoper
and homogeneous baselines (HET $-$ HOM) on the evaluation set. }
\label{tab:het_vs_hom}
\renewcommand{\arraystretch}{1.5}
\begin{tabular}{llccccccc}
\hline
& & & \multicolumn{6}{c}{Homogeneous baseline} \\
\cline{4-9}
Metric & Evaluation set & $n$ &
\multicolumn{2}{c}{RF} & \multicolumn{2}{c}{LR} & \multicolumn{2}{c}{MLP} \\
\cline{4-5}\cline{6-7}\cline{8-9}
& & & $\phi^{\mathrm{def}}=0.5$ & $\phi^{\mathrm{ROC}}_{r,m}$ & $\phi^{\mathrm{def}}=0.5$ & $\phi^{\mathrm{ROC}}_{r,m}$ & $\phi^{\mathrm{def}}=0.5$ & $\phi^{\mathrm{ROC}}_{r,m}$ \\
\hline
\multirow{4}{*}{Macro-$F_1$-score}
 & All    & 226 & 0.0047 & 0.0052 & 0.0301 & 0.1167 & 0.0207 & 0.0566 \\
 & Easy   & 74  & 0.0047 & 0.0179 & 0.0339 & 0.1181 & 0.0156 & 0.0394 \\
 & Medium & 24  & 0.0031 & 0.0128 & 0.0187 & 0.0670 & 0.0091 & 0.0649 \\
 & Hard   & 128 & 0.0052 & 0.0052 & 0.0249 & 0.1144 & 0.0208 & 0.0736 \\
\hline
\multirow{4}{*}{Macro-AP}
 & All    & 226 & \multicolumn{2}{c}{0.0038} & \multicolumn{2}{c}{0.0279}
 & \multicolumn{2}{c}{0.0209} \\
 & Easy   & 74  & \multicolumn{2}{c}{0.0126} & \multicolumn{2}{c}{0.0123}
 & \multicolumn{2}{c}{0.0204} \\
 & Medium & 24  & \multicolumn{2}{c}{0.0114} & \multicolumn{2}{c}{0.0129}
 & \multicolumn{2}{c}{0.0318} \\
 & Hard   & 128 & \multicolumn{2}{c}{0.0017} & \multicolumn{2}{c}{0.0377}
 & \multicolumn{2}{c}{0.0202} \\
\hline
\end{tabular}
\end{table}

\subsection{Downstream effects of learning-based reduction}\label{sec:results_downstream_learn}
This section examines the effect of FP predictions on the restricted rTPP, as the resulting route fixings may exclude routing alternatives needed to satisfy the operational constraints. We refer to an ML-induced route fixing arising from an FP prediction as an \emph{FP fixing}. \Cref{tab:fp_feasibility} shows that all $30$ instances without FP fixings remain feasible, whereas every instance requiring repair contains at least one FP fixing. However, repair incidence does not increase with the number of FPs: the repair share decreases from 57.45\% among instances with a single FP to 27.08\% among those with seven or more. The median fixing rate likewise exhibits no monotonic relationship with repair incidence. Hence, neither the number of FP fixings nor the overall fixing rate provides a reliable indication of whether repair is required.

Next, we examine the $77$ instances requiring fixing releases to quantify the extent to which FP fixings must be released to restore feasibility. \Cref{tab:fp_release} shows that, within this subset, instances with a single FP fixing require its release, whereas only a subset of FP fixings is generally released when multiple FPs are present. Across the $77$ instances, fewer than two FP fixings are released on average, accounting for 48.1\% of all FP fixings in this subset. Thus, restoring feasibility generally requires releasing only a limited subset of FP fixings, while the remaining ML-induced route fixings can be retained.

\begin{figure}[h]
    \centering
    \includegraphics[width=0.5\linewidth]{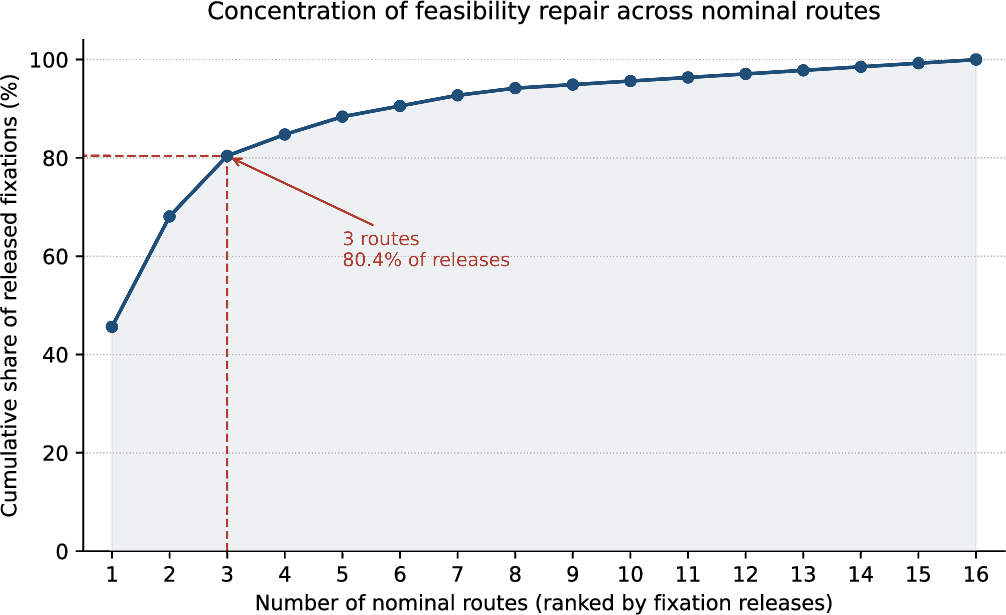}
    \caption{Cumulative concentration of released FP fixings across nominal route labels}
    \label{fig:cumulative_concentration_fp}
\end{figure}

The selective nature of these releases motivates examining their distribution across nominal route labels to determine whether they are broadly distributed or concentrated on particular routes.  \Cref{fig:cumulative_concentration_fp} shows that the three route labels most frequently associated with releases account for 80.4\% of all observed releases, while only 16 of the 72 nominal routes are associated with any release. Together, these results indicate that the feasibility impact of FP fixings depends more on which routing decisions are incorrectly fixed than on their number or overall fixing rate, with the required release concentrated on a small subset of nominal routes.


\begin{table}[h]
\centering
\caption{Relationship among FP fixings, fixing rate, and repair indices across the 226 instances.}
\label{tab:fp_feasibility}
\begin{tabular}{lcccc}
\hline
FP & Instances & Median fixing & Instances requiring & Repair \\
fixings & & rate (\%) & repair & share (\%) \\
\hline
0        & 30 & 62.74 &  0 &  0.00 \\
1        & 47 & 51.06 & 27 & 57.45 \\
2--3     & 53 & 59.75 & 20 & 37.74 \\
4--6     & 48 & 57.83 & 17 & 35.42 \\
$\geq7$ & 48 & 69.55 & 13 & 27.08 \\
\hline
All      & 226 & --- & 77 & 34.07 \\
\hline
\end{tabular}
\end{table}
\vspace{-1cm}
\begin{table}[h]
\centering
\caption{Extent of FP-fixing release by initial FP count across the instances requiring repair.}
\label{tab:fp_release}
\begin{tabular}{lccc}
\hline
FP          & Instances & Mean released & FP release \\
fixings   &           &  FP fixings & share  (\%) \\
\hline
$1$          & 27 & 1.00 & 100.0 \\
$2$--$4$     & 23 & 2.13 & 81.6 \\
$5$--$8$     & 19 & 1.95 & 33.0 \\
$\geq 9$     &  8 & 3.12 & 28.4 \\
\hline
All          & 77 & 1.79 & 48.1 \\
\hline
\end{tabular}
\end{table}

\subsection{Comparison of repair strategies}\label{sec:results_downstream_repair}
This section compares the conflict-guided repair, global repair, and CPLEX FeasOpt in terms of computational effort and preservation of the ML-induced reduction. We first examine repair time across disruption-duration groups and then assess pairwise repair-time performance across matched instances. Disruption duration provides a relevant stratification because longer capacity reductions affect a larger portion of the timetable and may increase the number of interacting train movements and conflicting routing decisions. Repair times are evaluated relative to the common $180$~s real-time budget.

\begin{figure}[t]
    \centering
    \includegraphics[height=8cm, width=0.75\linewidth]{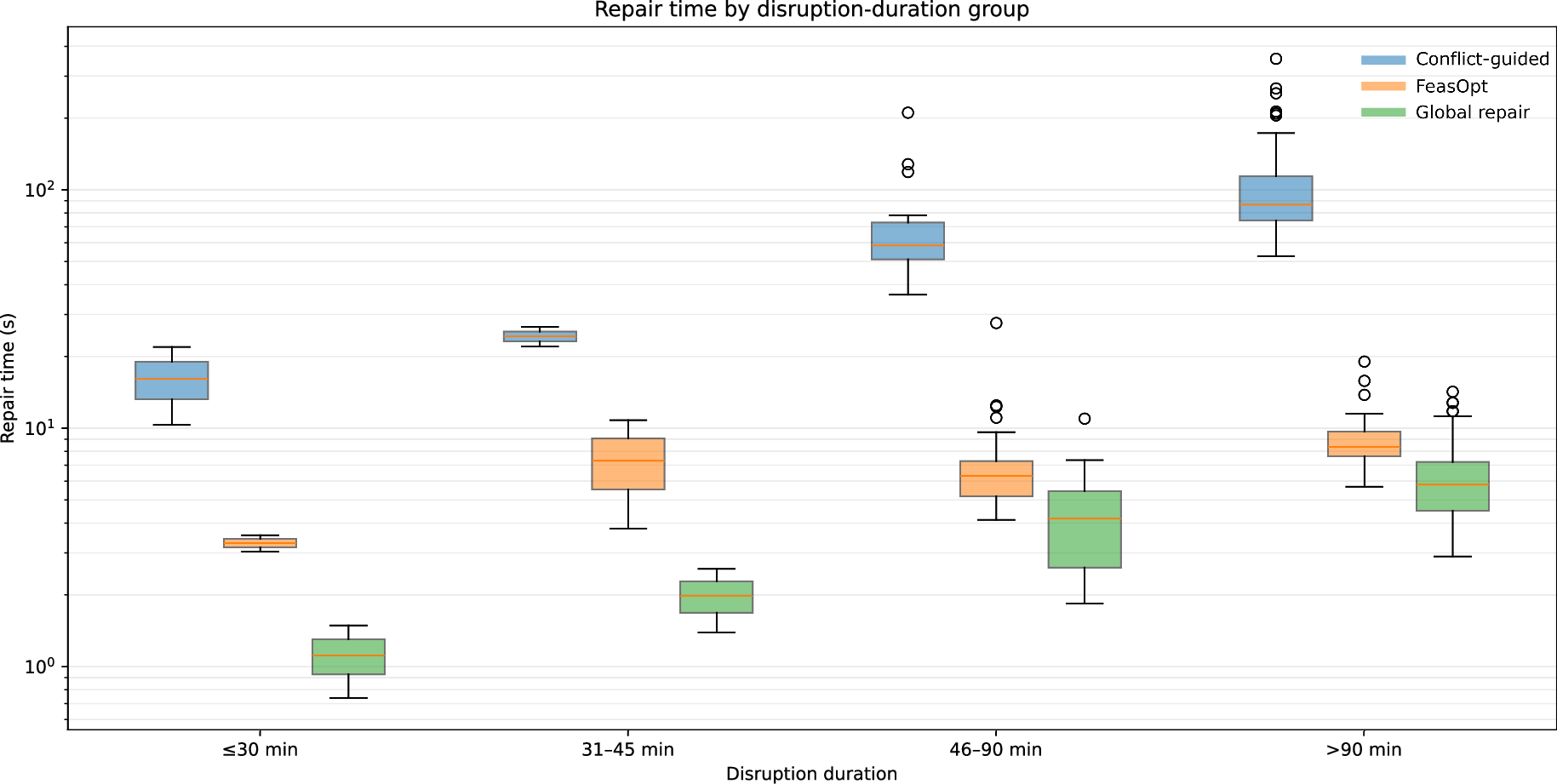}
    \caption{Computational effort of the different repair strategies across disruption-duration groups}
    \label{fig:repair_time_duration}
\end{figure}

\Cref{fig:repair_time_duration} shows that the differences in computational effort among the three repair strategies become more pronounced as disruption duration increases. For the 46--90~min and $>90$~min groups, conflict-guided repair exhibits substantially longer repair times and greater dispersion across instances, whereas global repair and FeasOpt remain comparatively stable. Further, timing decomposition shows that approximately 90\% of the conflict-guided repair time is spent in CPLEX conflict refinement, indicating that conflict identification, rather than solving the auxiliary repair MILP, is its main computational bottleneck. Global repair avoids repeated conflict identification by treating the ML-induced fixing set jointly under the cardinality restriction and exhibits only a moderate increase in repair time with disruption duration. The pairwise results in \Cref{tab:repair_pairwise} confirm that these differences are systematic: global repair and FeasOpt are faster than conflict-guided repair in all matched instances, while global repair is faster than FeasOpt in 81.9\% of cases. Overall, global repair provides the most consistently low repair times, particularly for longer disruptions.

\begin{table}[t]
\centering
\caption{Pairwise comparison of repair-time performance across matched instances.}
\label{tab:repair_pairwise}
\begin{tabular}{lcc}
\hline
Comparison & Faster strategy & Share of instances (\%) \\
\hline
Global repair vs.\ conflict-guided repair & Global repair & 100.0 \\
FeasOpt vs.\ conflict-guided repair       & FeasOpt       & 100.0 \\
Global repair vs.\ FeasOpt                & Global repair & 81.9 \\
                                           & FeasOpt       & 18.1 \\
\hline
\end{tabular}
\end{table}

We next compare the number of released fixings to assess how much of the ML-induced reduction is preserved after feasibility restoration. \Cref{tab:repair_release} shows that global repair and FeasOpt preserve the ML-induced reduction to essentially the same extent, consistent with both minimizing the number of fixing changes required to restore feasibility. Conflict-guided repair releases substantially more fixings because its minimum-release decisions are made sequentially within conflict-local subsets; releases selected in earlier iterations remain removed even if they would not be required in a globally coordinated repair. These local decisions can therefore accumulate across iterations, leading to a considerably larger relaxation. Together with the repair-time results, this shows that global repair reduces computational effort while preserving the ML-induced reduction. Accordingly, global repair is used in all subsequent 2RO and 3RO experiments.

\begin{table}[t]
\centering
\caption{Number of released ML-induced route fixings across repair strategies.}
\label{tab:repair_release}
\begin{tabular}{lcc}
\hline
Repair strategy & Mean releases & Median releases \\
\hline
Conflict-guided repair & 20.34 & 10 \\
Global repair          & 2.86 & 2 \\
FeasOpt                & 2.86 & 2 \\
\hline
\end{tabular}
\end{table}
\subsection{Effect of refinement on solution quality and primal progress}\label{sec:results_downstream_refine}
This section isolates the contribution of interaction-guided refinement by comparing 3RO with 2RO. Specifically, it examines whether selectively releasing additional route fixings improves primal solution quality sufficiently to offset the resulting expansion of the reduced routing search space within the fixed real-time budget. Results are reported both overall and conditional on the feasibility of the initially restricted rTPP, distinguishing refinement applied directly after reduction from refinement applied after feasibility restoration through repair.

This section isolates the contribution of interaction-guided refinement by comparing 3RO with 2RO. Specifically, it examines whether selectively releasing additional route fixings improves primal solution quality sufficiently to offset the resulting expansion of the restricted search space within the fixed real-time budget. Results are reported both overall and conditional on the feasibility of the initially restricted rTPP, distinguishing refinement applied directly after reduction from refinement applied after feasibility restoration through repair.

Performance is assessed along two complementary dimensions: terminal solution-quality improvement at $B=180$~s and primal progress, measured by the time required for 3RO to surpass the terminal 2RO incumbent. Let $z_{j'}^{\mathrm{2RO}}(\ell)$ and $z_{j'}^{\mathrm{3RO}}(\ell)$ denote the best incumbent objective values available at elapsed time, $\ell$ for instance $j'$. Terminal performance is evaluated at $\ell=B$, with 2RO as the reference; positive values indicate that 3RO attains a better terminal incumbent, whereas negative values indicate that refinement does not recover sufficient improvement within the available budget. To measure primal progress, let  $\pi_{j'}^{\mathrm{2RO}}=\min\{\ell: z_{j'}^{\mathrm{2RO}}(\ell)\leq z_{j'}^{\mathrm{2RO}}(B)\}$ denote the first time at which 2RO attains its terminal incumbent, and let $\ell_{j'}^{\mathrm{ref}}$ denote the start time of refinement. We define $\pi_{j'}^{\mathrm{3RO},+}=\min\{\ell\in[ \ell_{j'}^{\mathrm{ref}}, B]: z_{j'}^{\mathrm{3RO}}(\ell)< z_{j'}^{\mathrm{2RO}}(B)\}$ as the first post-refinement time at which 3RO strictly improves upon the terminal 2RO incumbent. The attainment-time difference is defined as $\Delta L_{j'}=\pi_{j'}^{\mathrm{3RO},+}-\pi_{j'}^{\mathrm{2RO}}$. A negative value indicates that 3RO surpasses the benchmark before 2RO attains it, whereas a positive value indicates the additional time required by 3RO. If 3RO does not surpass the benchmark within $B$, $\Delta L_{j'}$ is right-censored at $B-\pi_{j'}^{\mathrm{2RO}}$; that is, $\Delta L_{j'}> B-\pi_{j'}^{\mathrm{2RO}}$.

\Cref{tab:refinement_impact} shows that refinement contributes to both feasibility groups. When the learned restriction is initially feasible, it improves the terminal incumbent in 75.6\% of instances, with a median improvement of 38.69\%, indicating that feasibility alone does not ensure that the restricted MILP retains sufficient routing alternatives for high-quality solutions. After repair, refinement still improves 66.0\% of instances, with only one deterioration. These results separate the two adaptive roles of 3RO: repair restores feasibility, whereas refinement recovers additional routing alternatives for solution-quality improvement.

\begin{table}[t]
\centering
\caption{Terminal effect of refinement by initial reduced-model feasibility}
\label{tab:refinement_impact}
\begin{tabular}{lrrrrr}
\hline
& & \multicolumn{3}{c}{3RO relative to 2RO} & \\
\cline{3-5}
Initial reduced model & $n$ & Better & Equal & Worse &
Median improvement (\%)\\
\hline
Feasible        & 78  & 59 & 9  & 10 & 38.69\\
Requires repair & 50  & 33 & 16 & 1  & 31.11\\
\hline
All             & 128 & 92 & 25 & 11 & 34.00\\
\hline
\end{tabular}
\end{table}

\begin{figure}[t]
    \centering
    \includegraphics[height=10cm, width=0.85\linewidth]{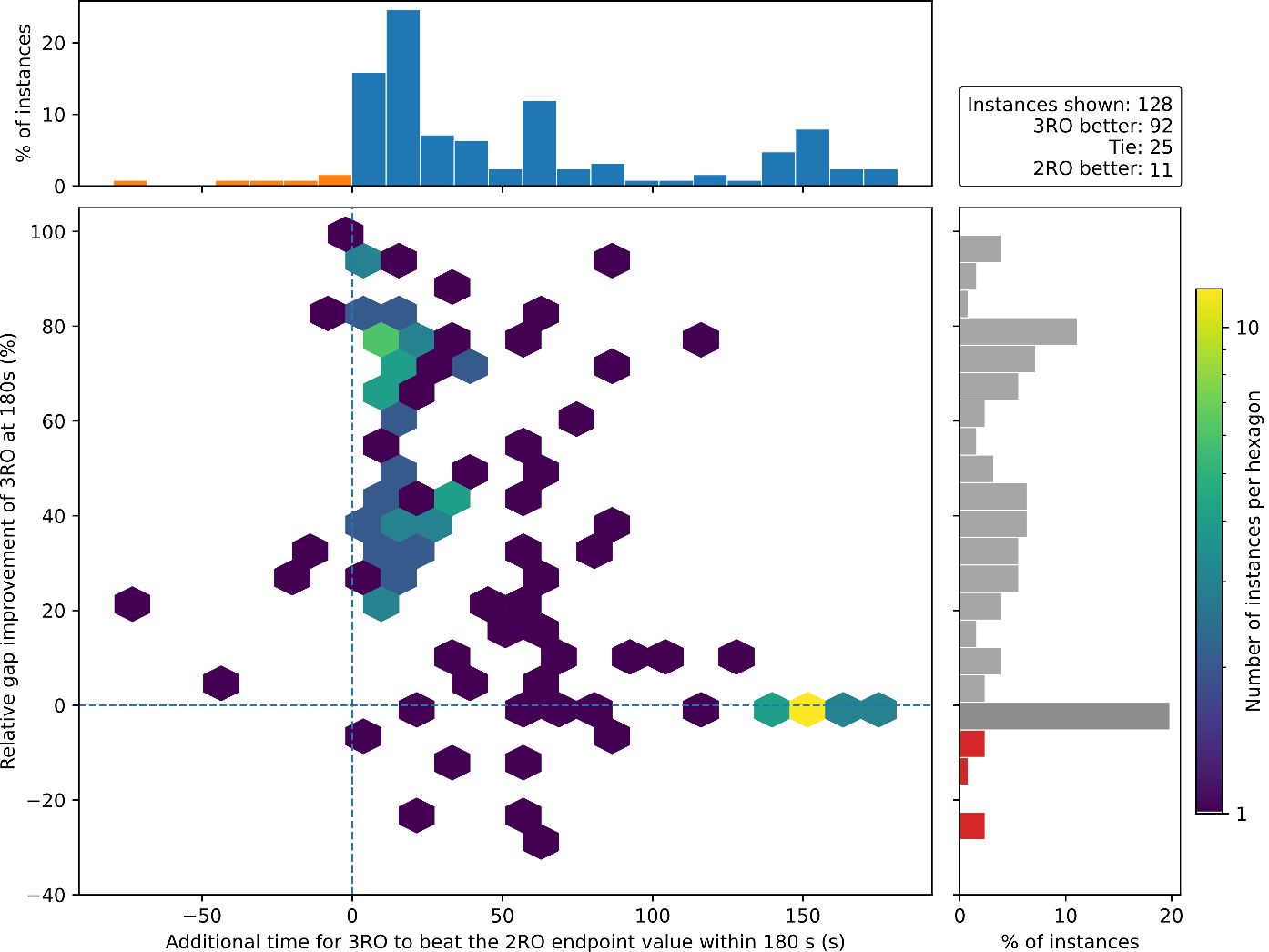}
    \caption{Terminal solution-quality improvement and attainment-time difference of 3RO relative to 2RO}
    \label{fig:endpoint_improvement_3RO_2RO}
\end{figure}
Having established the terminal effect of refinement across the two feasibility groups, \Cref{fig:endpoint_improvement_3RO_2RO} examines how the magnitude of these improvements relates to the time required to attain them across the 128 hard instances. Among instances for which 3RO surpasses the terminal 2RO incumbent, 55\% attain this improvement with an attainment-time difference between 0 and 50~s. Six instances exhibit $\Delta L_{j'}<0$, indicating that 3RO surpasses the terminal 2RO incumbent before 2RO itself first attains that value. Some improvements nevertheless require substantially longer computation. Instances with zero or negative terminal improvement are right-censored with respect to the attainment-time measure because 3RO does not surpass the terminal 2RO benchmark within $B$. These results illustrate the trade-off introduced by refinement: releasing additional route fixings can improve the attainable solution quality, but the expansion of the reduced routing search space must be exploited within the remaining real-time budget. 

\begin{figure}[t]
    \centering
    \includegraphics[width=0.9\linewidth]{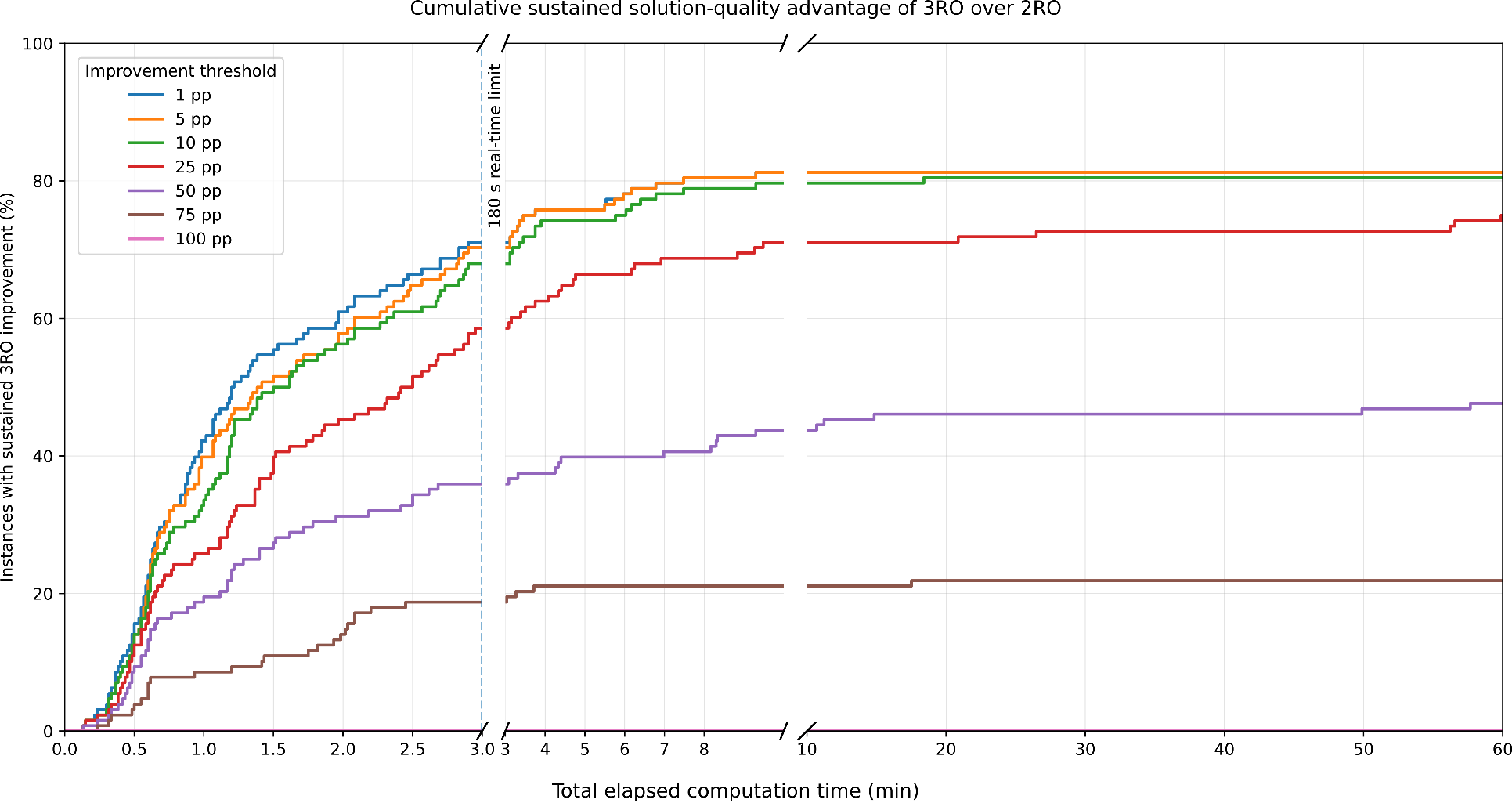}
    \caption{Cumulative sustained solution-quality advantage of 3RO over 2RO}
    \label{fig:sustained_advantage_2ro_3ro}
\end{figure}

To assess whether the refinement gains observed within the 180~s real-time horizon persist with additional computation, we extend the comparison diagnostically to 60~min.  \Cref{fig:sustained_advantage_2ro_3ro} reports the evolution of the primal advantage of 3RO over 2RO over this extended horizon. For each improvement threshold $\theta$, the curve gives the share of instances for which 3RO attains an advantage of at least $\theta$ percentage points and maintains this advantage thereafter. The resulting curves are therefore monotone and capture sustained rather than transient improvements. Within the 180~s real-time budget, refinement already provides a sustained advantage over 2RO for a large share of the instances: approximately 71\% achieve at least 1~pp improvement, 59\% at least 25~pp, and 36\% at least 50~pp. Over the extended 60-min horizon, these shares increase to approximately 81\%, 74\%, and 48\%, respectively. Thus, the gains from refinement are frequently substantial and persist beyond the real-time horizon. 

Most of the increase in sustained advantage occurs within the first several minutes. By approximately 10~min, the 1, 5, 10, and 25~pp curves are already close to their 60-min levels, while the remaining curves show only modest additional growth thereafter. This indicates diminishing returns from extending the computation horizon, although some instances still realize larger gains beyond the real-time budget. The 1~pp curve stabilizes at approximately 81\% because, for the remaining 19\% of instances, 2RO already attains the unrestricted-MILP optimum. This indicates that the ML-induced fixings preserve the routing decisions required for an optimal solution in these cases, leaving no further objective improvement available to refinement.

\subsection{Comparison with heuristic and unrestricted optimization benchmarks}\label{sec:results_sol_quality}
To position the computational performance of 3RO, we compare it with two complementary benchmarks: the recent rTPP heuristic of \citet{garcia2024efficient} and the unrestricted MILP. The following comparisons examine the time-to-quality performance of 3RO relative to these two alternatives.

\noindent\textbf{Comparison with a recent rTPP heuristic.}
The greedy heuristic emphasizes rapid feasibility attainment, requiring on average 0.7~s to obtain a feasible solution and less than 2~s to reach its final solution. This speed, however, does not translate into complete feasibility coverage: the heuristic returns a feasible solution for 216 of the 226 evaluation instances, whereas 3RO obtains feasible solutions for all instances within the real-time budget. Although 3RO does not match this immediate response, it closes the feasibility gap rapidly, with a median time to first feasible incumbent of 4.6~s. As shown in \Cref{tab:garcia_quality_over_time}, about 90\% of comparable instances have a feasible 3RO incumbent within 15~s, and all do so within one minute. Solution-quality gains emerge over the same early period: by 15~s, 3RO matches or improves upon the heuristic in more than half of the instances, with about one quarter achieving improvements of at least 25\%. 

\begin{table}[h]
\centering
\caption{Evolution of 3RO feasibility and solution-quality relative to the greedy heuristic.}
\label{tab:garcia_quality_over_time}
\begin{tabular}{rrrrrrrrr}
\hline
Time & Feasible 3RO & Better & Equal & \multicolumn{3}{c}{Improvement}& \multicolumn{2}{c}{Loss} \ \\
(s)  & (\%) & (\%) & (\%) & $\geq25\%$ & $\geq50\%$ & $\geq75\%$& $\geq25\%$ & $\geq50\%$ \\
\hline
5   & 52.3  & 11.7 &20.6& 11.7 &  7.5 &  4.7& --- & --- \\
15  & 91.1  & 29.0 &26.6& 25.7 & 17.8 & 10.3& --- & --- \\
30  & 96.3  & 45.3 &27.6& 40.2 & 30.8 & 21.0& --- & --- \\
60  & 100.0 & 53.7 &27.1& 46.7 & 36.0 & 28.0& 15.5 & 13.2 \\
120 & 100.0 & 59.3 &26.2& 56.1 & 43.5 & 31.3& 11.5 & 9.3 \\
180 & 100.0 & 63.1 &25.7& 58.9 & 46.3 & 33.6& 8.4 & 6.2 \\
\hline
\end{tabular}
\end{table}

The solution-quality profile shifts increasingly in favor of 3RO within the first minute. By 60~s, 3RO matches or improves on the heuristic in more than 80\% of instances. More importantly, the gains are not concentrated around marginal improvements: among the instances in which 3RO outperforms the heuristic, approximately two-thirds achieve an improvement of at least 50\%, and more than half improve by at least 75\%. Thus, when 3RO improves upon the heuristic, the gain is frequently substantial rather than marginal. The inferior tail contracts as additional computation becomes available, indicating that some early quality disadvantages are resolved as the search progresses. Loss magnitudes are reported from 60~s onward because, before this point, some instances do not yet have a feasible 3RO incumbent and therefore do not admit a relative solution-quality comparison. Overall, the two approaches exhibit distinct computational profiles: the heuristic prioritizes near-immediate feasibility when successful, whereas 3RO requires more computation but achieves complete feasibility coverage and increasingly favors higher-quality incumbents within the real-time budget.

\noindent\textbf{Comparison with the unrestricted MILP.}
We next compare 3RO with the unrestricted MILP to assess the computational effect of restricting the search space. The comparison is conducted across all 226 instances and, separately, the 128 hard instances, considering two computation horizons. The common 180~s horizon compares both approaches under equal wall-clock effort, while the 3,600~s benchmark quantifies the solution-quality trade-off relative to unrestricted optimization when substantially more computation time is available.

\begin{figure}[h]
    \centering
    \includegraphics[height= 10cm, width=1.0\linewidth]{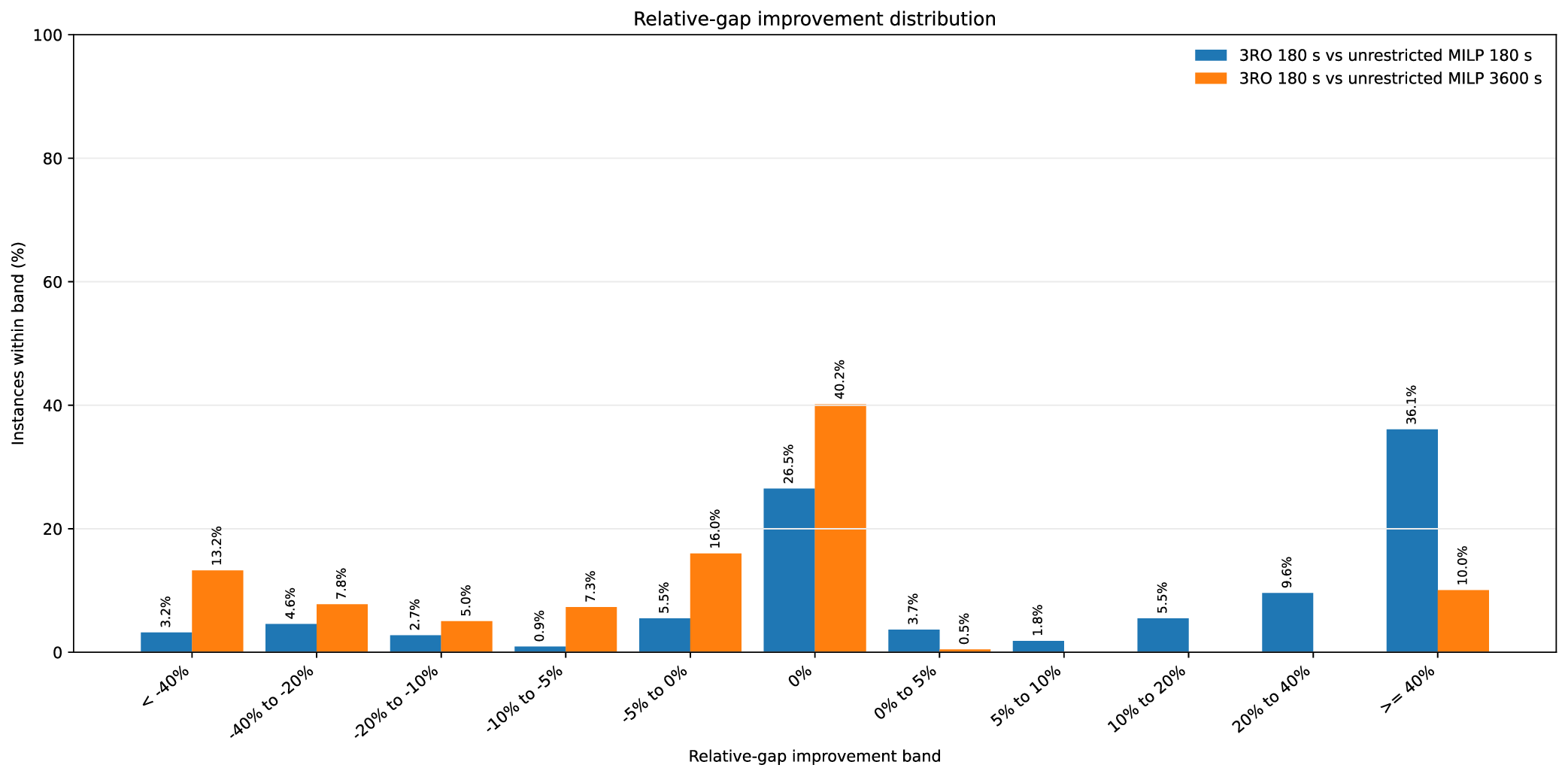}
    \caption{Relative solution-quality improvement of 3RO over full-MILP benchmarks: All instances}
    \label{fig:all_instances_relative_imp}
\end{figure}

\begin{figure}[h]
    \centering
    \includegraphics[height= 10cm, width=1.0\linewidth]{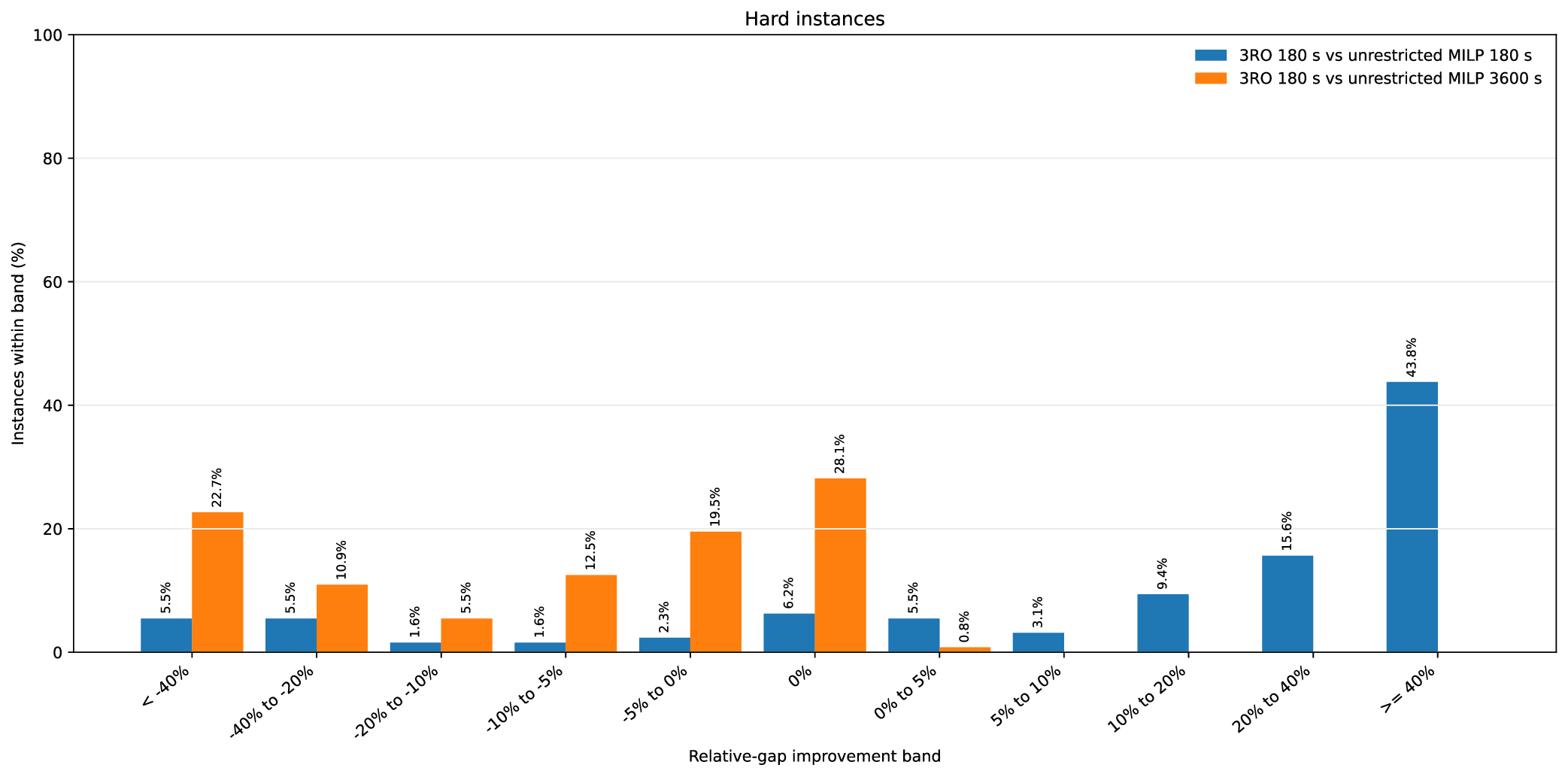}
    \caption{Relative solution-quality improvement of 3RO over full-MILP benchmarks: Hard instances}
    \label{fig:hard_instances_relative_imp}
\end{figure}

As shown in \Cref{fig:all_instances_relative_imp,fig:hard_instances_relative_imp}, at the common 180~s horizon, 3RO obtains a strictly better incumbent than the unrestricted MILP in 56.7\% of all instances and 77.4\% of the hard instances. The gains are substantial in many cases: improvements of at least 40\% occur in 36.1\% of all instances and 43.8\% of the hard instances. Including ties, 3RO is no worse than the unrestricted MILP in approximately 83\% of both groups. The larger share of improvements among the hard instances indicates that, within the real-time computation horizon, the benefit of search-space restriction becomes more pronounced as instance difficulty increases. 

The same comparison shows that the frequency of inferior 3RO solutions, indicated by negative relative-gap, remains nearly unchanged, at 16.9\% across all instances and 16.5\% among hard instances. The severity of these losses, however, increases with instance difficulty: losses exceeding 20\% occur in 7.8\% of all instances and 11.0\% of the hard instances. Thus, the restricted search is not uniformly beneficial; when relevant routing alternatives remain excluded, the resulting solution-quality loss can be more pronounced for computationally difficult instances. Overall, instance difficulty amplifies the effect of search-space restriction: the frequency of deterioration remains stable, while both the observed gains and the magnitude of occasional losses increase.

Next, when the unrestricted MILP is allowed the extended 3,600~s horizon, the advantage of 3RO narrows, particularly for the hard instances, as shown in \Cref{fig:all_instances_relative_imp,fig:hard_instances_relative_imp}. Nevertheless, the 180~s 3RO incumbent already matches the 3,600~s full-MILP incumbent in 28.1\% of the hard instances, with a further 19.5\% remaining within 5\%. The stronger negative tail under the extended benchmark indicates that unrestricted optimization can recover additional solution quality when given substantially more computation time. Taken together, the two horizons characterize the intended trade-off of 3RO: restricting the search improves access to high-quality incumbents within the real-time budget, while potentially foregoing solutions that become attainable through prolonged unrestricted search.

\subsection{Conclusion and future directions}
This study develops 3RO, a reduce--repair--refine--optimize framework for real-time train platforming under disruptions. The framework uses offline learning to identify route fixings that reduce the routing search space while retaining the original rTPP objective and operational constraints for determining the rescheduling decisions. When the initial reduction renders the restricted rTPP infeasible, repair selectively releases route fixings to restore feasibility. Subsequently, interaction-guided refinement releases additional route fixings associated with operationally relevant trains to recover routing alternatives that may improve solution quality. Thus, 3RO adapts the reduced routing search space during the online solution process rather than relying solely on the initial learning-based reduction.

The computational results confirm that predictive performance alone does not fully characterize the effectiveness of learning-based reduction. Although the heterogeneous Scoper consistently outperforms the homogeneous variants, the downstream impact of FP fixings is strongly route-dependent: only 1.79 FP fixings are released on average to restore feasibility, and three route labels account for 80.4\% of all observed releases. Interaction-guided refinement further improves the terminal incumbent in 71.9\% of the hard instances, indicating that selectively recovering routing alternatives generally outweighs the additional computational effort associated with expanding the reduced routing search space.  Relative to the recent greedy rTPP heuristic, 3RO requires more time to establish feasibility but closes this gap rapidly and, within the first minute, matches or improves upon the heuristic in more than 80\% of the instances, with substantial improvements in many cases. Finally, against the unrestricted MILP under the real-time budget, the asymmetric distribution of incumbent improvement reveals a generally favorable but non-uniform time-to-quality trade-off: reducing the routing search space frequently directs computational effort toward high-quality incumbents, while the relatively small negative tail reflects the residual risk of excluding improving routing alternatives. Although this advantage narrows against the 3,600~s unrestricted MILP benchmark, the 180~s 3RO incumbent remains within $\pm5\%$ in 47.6\% of the hard instances. Thus, nearly half of hard instances attain comparable solution quality using one-twentieth of the computation time, while the remaining negative tail reflects the solution-quality cost of search-space reduction over longer horizons. 

The findings also highlight a few limitations. First, route-stability labels are derived from a single reference solution and may therefore classify routes used in alternative optimal or near-optimal solutions as unstable. Second, the observed release patterns are conditional on the selected repair strategy: alternative minimum-cardinality release sets may retain different route fixings and consequently produce different post-repair feasible regions and downstream solution quality. Finally, the 3,600~s unrestricted MILP solution remains a finite-time incumbent; the extended comparison therefore assesses time-to-quality rather than proximity to global optimality. 

These limitations suggest several directions for further research. Route stability could be derived from multiple reference solutions or solution pools to better represent alternative high-quality routings and support probabilistic or confidence-based labels. Extending the Scoper toward structured or graph-based models could capture dependencies among trains, routes, and infrastructure, while calibration and decision-aware training objectives could account for the asymmetric downstream impact of FP fixings. During the online phase, lexicographic repair could use the rTPP objective to break ties among minimum-cardinality release sets, thereby favoring repairs with greater downstream solution potential without increasing the number of releases. Similarly, refinement could incorporate adaptive release budgets and optimization- or learning-based candidate selection to target route fixings with high improvement potential. Finally, finite-time incumbent comparisons could be complemented by solver bounds, optimality gaps, and proven-optimal benchmarks for tractable instances.

\bibliographystyle{authordate1} 
\bibliography{references}
\setlength{\bibsep}{0pt}

\clearpage

\appendix

\section{MILP formulation of the rTPP}
\label{app:rtpp_formulation}

This appendix presents the complete MILP formulation of the rTPP underlying the 3RO framework and is used as the benchmark in the computational experiments. The sets, parameters, and decision variables are summarized in Tables~\ref{tab:notation_pars} and~\ref{tab:notation_vars}. 

\begin{table}[h]
\captionsetup{singlelinecheck=false, justification=raggedright}
\caption{Notation of sets and parameters used in the MILP.}
\label{tab:notation_pars}
\footnotesize
\renewcommand{\arraystretch}{1.5}
\begin{tabular}{@{}lp{0.75\textwidth}@{}}
\hline
\textbf{Notation} & \textbf{Description} \\
\hline

\multicolumn{2}{@{}l}{\textbf{Sets}} \\
$\mathcal{T}$ & Scheduled trains \\
$\mathcal{T}_o, \mathcal{T}_s$ & Originating and terminating trains \\
$\widetilde{\mathcal{T}}, \mathcal{\bar{T}}$ & Train pairs running in the same / opposing direction that may generate conflicts \\
$\mathcal{Q}$ & Running directions of trains \\
$\mathcal{R}_t^{+}, \mathcal{R}_t^{-}$ & Feasible inbound and outbound routes for train $t$ \\
$\mathcal{R}^x$ & Dummy routes used for cancellation \\
$\mathcal{P}, \mathcal{P}_t$ & Platforms in the station area, feasible platforms for train $t$ \\

\hline
\multicolumn{2}{@{}l}{\textbf{Parameters}} \\
$\lambda_X, \lambda_D, \lambda_P, \lambda_R \in \mathbb{R}_{\ge 0}$ & Weights of the cost components in the objective function \\
$q_t$ & Running direction of train $t$ \\
$\tau^{+}_t, \tau^{-}_t$ & Inbound and outbound running time of train $t$ \\
$p(r)$ & Platform associated with route $r$ \\
$\mu_{r\dot{r}}$ & 1 if inbound routes $r$ and $\dot{r}$ share the same entry vertex $e_h$; 0 otherwise \\
$\nu_{r\dot{r}}$ & 1 if $\mathcal{V}(r) \cap \mathcal{V}(\dot{r}) \setminus \mathcal{V}^{\mathrm{b}} \neq \emptyset$ for inbound routes $r$ and $\dot{r}$; 0 otherwise \\
$\pi_{r\dot{r}}$ & 1 if inbound routes $r$ and $\dot{r}$ terminate at the same platform track; 0 otherwise \\
$\omega_{r\dot{r}}$ & 1 if $\mathcal{V}(r) \cap \mathcal{V}(\dot{r}) \neq \emptyset$ for outbound routes $r$ and $\dot{r}$; 0 otherwise \\
$\zeta_{r\dot{r}}$ & 1 if $\mathcal{V}(r) \cap \mathcal{V}(\dot{r}) \neq \emptyset$ for opposing routes $r$ and $\dot{r}$; 0 otherwise \\
$\eta_{r\dot{r}}$ & 1 if opposing routes $r$ and $\dot{r}$ intersect only at a platform track; 0 otherwise \\
$\hat{\alpha}^r_t$ & 1 if train $t$ is nominally assigned to inbound route $r$; 0 otherwise \\
$\hat{\beta}^r_t$ & 1 if train $t$ is nominally assigned to outbound route $r$; 0 otherwise \\
$\hat{p}_t$ & Nominally assigned platform of train $t$ \\
$\hat{a}^e_t$, $\hat{d}^e_t$ & Nominal arrival and departure times at the boundary vertex ($\mathcal{V}^{\mathrm{b}}$) \\
$\hat{a}^p_t$, $\hat{d}^p_t$ & Nominal arrival and departure times at the platform track \\
$\hat{w}^e_t, \hat{w}^p_t$ & Nominal waiting time at the entry vertex and dwell time at the platform track \\
$\bar{w}^e, \bar{w}^{p}$ & Maximum allowed waiting time at the entry vertex and dwell time at the platform track \\
$M$ & Sufficiently large constant \\
\hline
\end{tabular}
\end{table}

\begin{table}[H]
\captionsetup{singlelinecheck=false, justification=raggedright}
\caption{Notation of decision variables in the MILP.}
\label{tab:notation_vars}
\footnotesize
\renewcommand{\arraystretch}{1.5}
\begin{tabular}{@{}lp{0.75\textwidth}@{}}
\hline
\textbf{Notation} & \textbf{Description} \\
\hline

\multicolumn{2}{@{}l}{\textbf{Binary decision variables}} \\
$\alpha^{r}_{t}$ & 1 if train $t$ is assigned to inbound route $r$; 0 otherwise \\
$\beta^{r}_{t}$ & 1 if train $t$ is assigned to outbound route $r$; 0 otherwise \\
$\gamma^{r}_{t}$ & 1 if train $t$ is canceled i.e. assigned dummy route $r \in \mathcal{R}^{x}$; 0 otherwise \\
$\delta^{+}_t, \delta^{-}_t$ & Indicator for inbound / outbound route change of train $t$ \\
$y_{t\dot{t}}$ & 1 if train $t$ enters the station area before train $\dot{t}$ at a shared entry vertex; 0 otherwise \\
$z_{t\dot{t}}$ & 1 if train $t$ departs before train $\dot{t}$ from a shared platform track; 0 otherwise \\

\hline
\multicolumn{2}{@{}l}{\textbf{Continuous decision variables}} \\
$a^e_t, d^e_t$ & Revised arrival and departure time of train $t$ at the boundary vertex  \\
$a^p_t, d^p_t$ & Revised arrival and departure time of train $t$ at the platform track \\
$\Delta p_t$ & Platform deviation of train $t$ from its nominal assignment \\
$\Delta w^e_t, \Delta w^p_t$ & Deviation in waiting time at the entry vertex and in dwell time at the platform track \\

\hline
\end{tabular}
\end{table}

\noindent \textbf{Objective function.}
\begin{equation}
\begin{aligned}
\min \quad & \lambda_X \sum_{t \in \mathcal{T}} \sum_{r \in \mathcal{R}^x} \gamma^r_t
+ \lambda_D \sum_{t \in \mathcal{T}} \left( \Delta w^e_t + \Delta w^p_t \right) \\
& + \lambda_P \sum_{t \in \mathcal{T}} \Delta p_t
+ \lambda_\mathcal{R} \sum_{t \in \mathcal{T}} \left( \delta^{+}_t + \delta^{-}_t \right),
\end{aligned}
\label{eq:full_model_obj}
\end{equation}

The objective function minimizes the total weighted operational deviation from the nominal timetable. Here, $\gamma^r_t$ indicates cancellation of train $t$ via assignment to dummy route $r\in\mathcal{R}^x$; $\Delta w_t^e$ and $\Delta w_t^p$ represent additional waiting time at the entry vertex and additional dwell time at the platform track, respectively; $\Delta p_t$ measures the deviation from the nominal platform-track assignment; and $\delta_t^+$ and $\delta_t^-$ indicate changes from the nominal inbound and outbound routes, respectively. The weights represent the relative importance of cancellations ($\lambda_X$), delay-related deviations ($\lambda_D$), platform deviations ($\lambda_P$), and route changes ($\lambda_{\mathcal{R}}$), with $\lambda_X \gg \lambda_D > \lambda_P > \lambda_{\mathcal{R}}$. Timing constraints (\ref{eq:full_model_2})--(\ref{eq:full_model_6}) enforce consistency among arrival, departure, dwell, and waiting times. Flow conservation constraints (\ref{eq:full_model_7})--(\ref{eq:full_model_9}) ensure that each train is either cancelled or assigned a compatible pair of inbound and outbound routes. Headway constraints (\ref{eq:full_model_10})--(\ref{eq:full_model_20}) enforce the required temporal separation between conflicting train movements over shared infrastructure. Operational deviation constraints (\ref{eq:full_model_21})--(\ref{eq:full_model_26}) quantify changes from the nominal timetable in routing, platforming, and timing decisions. Finally, domain definitions (\ref{eq:full_model_27})--(\ref{eq:full_model_30}) specify the integrality and feasible bounds of the binary and continuous decision variables.

\noindent \textbf{Arrival and departure time constraints.}
\begin{align}
&a^e_t \geq \hat{a}^e_t - M\sum_{r\in \mathcal{R}^0}\gamma^r_t &\forall t\in \mathcal{T}\label{eq:full_model_2}\\
&a^p_t \geq \hat{ap}_t - M\sum_{r\in \mathcal{R}^0}\gamma^r_t &\forall t\in \mathcal{T}\label{eq:full_model_3}\\
&d^e_t = a^e_t + \hat{w}^e_t + \Delta w^e_t   &\forall t\in \mathcal{T}\label{eq:full_model_4}\\
&a^p_t = d^e_t + \tau^{+}_t &\forall t\in \mathcal{T}\label{eq:full_model_5}\\
&d^p_t = a^p_t + \hat{w}^p_t + \Delta w^p_t  &\forall t\in \mathcal{T}\label{eq:full_model_6} 
\end{align}

\noindent \textbf{Flow conservation.}
\begin{align}
&\sum_{r\in R^{+}_t|p(r)=p} \alpha^{r}_t - \sum_{r\in R^{-}_t|p(r)=p} \beta^{r}_t = 0 &\forall t\in \mathcal{T}, p\in \mathcal{P}_t\label{eq:full_model_7}\\
&\sum_{r\in \mathcal{R}^{+}_t} \alpha^{r}_t + \sum_{r\in \mathcal{R}^{0}} \gamma^{r}_t = 1 &\forall t\in \mathcal{T}_{o}\label{eq:full_model_8}\\
&\sum_{r\in \mathcal{R}^{+}_t} \alpha^{r}_t = 1 &\forall t\in \mathcal{T}\setminus \mathcal{T}_{o}\label{eq:full_model_9}
\end{align}

\noindent\textbf{Inbound headway constraints.}
\begin{align}
& d^e_{\dot{t}} \geq d^e_{t} + \epsilon - M[2 - \alpha^r_t - \alpha^{\dot{r}}_{t}] \notag \\
& \qquad \forall (t,\dot{t}) \in \widetilde{T}_h,\; r \in \mathcal{R}^{+}_t,\; \dot{r} \in \mathcal{R}^{+}_{\dot{t}},\; \mu_{r\dot{r}} = 1
\label{eq:full_model_10} \\[2pt]
& d^e_{\dot{t}} \geq a^p_{t} + \epsilon - M[3 - y_{t\dot{t}} - \alpha^r_t - \alpha^{\dot{r}}_{\dot{t}}] \notag \\
& \qquad \forall (t,\dot{t}) \in \widetilde{T},\; r \in \mathcal{R}^{+}_t,\; \dot{r} \in \mathcal{R}^{+}_{\dot{t}},\; \nu_{r\dot{r}} = 1
\label{eq:full_model_11} \\[2pt]
& d^e_{t} \geq a^p_{\dot{t}} + \epsilon - M[2 + y_{t\dot{t}} - \alpha^r_t - \alpha^{\dot{r}}_t] \notag \\
& \qquad \forall (t,\dot{t}) \in \widetilde{T},\; r \in \mathcal{R}^{+}_t,\; \dot{r} \in \mathcal{R}^{+}_{\dot{t}},\; \nu_{r\dot{r}} = 1
\label{eq:full_model_12} \\[2pt]
& d^e_{\dot{t}} \geq d^p_{t} + \epsilon - M[3 - y_{t\dot{t}} - \alpha^{r}_t - \alpha^{\dot{r}}_{\dot{t}}] \notag \\
& \qquad \forall (t,\dot{t}) \in \widetilde{T},\; r \in \mathcal{R}^{+}_t,\; \dot{r} \in \mathcal{R}^{+}_{\dot{t}},\; \pi_{r\dot{r}} = 1
\label{eq:full_model_13} \\[2pt]
& d^e_{t} \geq d^p_{\dot{t}} + \epsilon - M[2 + y_{t\dot{t}} - \alpha^r_t - \alpha^{\dot{r}}_{\dot{t}}] \notag \\
& \qquad \forall (t,\dot{t}) \in \widetilde{T},\; r \in \mathcal{R}^{+}_t,\; \dot{r} \in \mathcal{R}^{+}_{\dot{t}},\; \pi_{r\dot{r}} = 1
\label{eq:full_model_14}
\end{align}

\noindent\textbf{Outbound headway constraints.}
\begin{align}
& d^p_{t} \geq d^p_{\dot{t}} + \tau^{-}_{\dot{t}} + \epsilon - M[3 - z_{t\dot{t}} - \beta^r_t - \beta^{\dot{r}}_{\dot{t}}] \notag \\
& \qquad \forall (t,\dot{t}) \in \widetilde{T},\; r \in \mathcal{R}^{-}_t,\; \dot{r} \in \mathcal{R}^{-}_{\dot{t}},\; \omega_{r\dot{r}} = 1
\label{eq:full_model_15} \\[2pt]
& d^p_{\dot{t}} \geq d^p_{t} + \tau^{-}_{t} + \epsilon - M[2 + z_{t\dot{t}} - \beta^r_t - \beta^{\dot{r}}_{\dot{t}}] \notag \\
& \qquad \forall (t,\dot{t}) \in \widetilde{T},\; r \in \mathcal{R}^{-}_t,\; \dot{r} \in \mathcal{R}^{-}_{\dot{t}},\; \omega_{r\dot{r}} = 1
\label{eq:full_model_16}
\end{align}

\noindent\textbf{Crossover headway constraints.}
\begin{align}
& d^e_{t} \geq d^p_{\dot{t}} + \tau^{-}_{\dot{t}} + \epsilon - M[3 - z_{t\dot{t}} - \alpha^{\dot{r}}_{\dot{t}} - \beta^r_t] \notag \\
& \qquad \forall (t,\dot{t}) \in \overline{T},\; \dot{r} \in \mathcal{R}^{+}_{\dot{t}},\; r \in \mathcal{R}^{-}_t,\; \zeta_{r\dot{r}} = 1
\label{eq:full_model_17} \\[2pt]
& d^e_{\dot{t}} \geq d^p_{t} + \tau^{-}_{t} + \epsilon - M[2 + z_{t\dot{t}} - \alpha^{r}_t - \beta^{\dot{r}}_{\dot{t}}] \notag \\
& \qquad \forall (t,\dot{t}) \in \overline{T},\; r \in \mathcal{R}^{+}_t,\; \dot{r} \in \mathcal{R}^{-}_{\dot{t}},\; \zeta_{r\dot{r}} = 1
\label{eq:full_model_18} \\[2pt]
& d^e_{t} \geq d^p_{\dot{t}} + \epsilon - M[3 - z_{t\dot{t}} - \alpha^{r}_{t} - \beta^{\dot{r}}_{\dot{t}}] \notag \\
& \qquad \forall (t,\dot{t}) \in \overline{T},\; r \in \mathcal{R}^{+}_t,\; \dot{r} \in \mathcal{R}^{-}_{\dot{t}},\; \eta_{r\dot{r}} = 1
\label{eq:full_model_19} \\[2pt]
& d^e_{\dot{t}} \geq d^p_{t} + \epsilon - M[2 + z_{t\dot{t}} - \alpha^{\dot{r}}_{\dot{t}} - \beta^{r}_{t}] \notag \\
& \qquad \forall (t,\dot{t}) \in \overline{T},\; r \in \mathcal{R}^{-}_t,\; \dot{r} \in \mathcal{R}^{+}_{\dot{t}},\; \eta_{r\dot{r}} = 1
\label{eq:full_model_20}
\end{align}

\noindent \textbf{Operational deviation constraints.}
\begin{align}
&\Delta{p}_t \geq \alpha^t_r \times |p(r) - p^0_t| &\forall t\in \mathcal{T},r\in \mathcal{R}^+_t\label{eq:full_model_21}\\
&|\alpha^r_t - \hat{\alpha}^r_t| \leq M \delta^+_t  &\forall t\in \mathcal{T}, r\in \mathcal{R}^+_t\label{eq:full_model_22}\\
&|\beta^r_t - \hat{\beta}^r_t| \leq M \delta^-_t  &\forall t\in \mathcal{T}, r\in \mathcal{R}^-_t\label{eq:full_model_23}
\end{align}

\noindent \textbf{Bounds on operational deviation.}
\begin{align}
&\Delta w^e_t \leq \bar{w}^e - \hat{w}^e_t  &\forall t\in \mathcal{T}\label{eq:full_model_24}\\
&\Delta w^p_t \leq \bar{w}^p - \hat{w}^p_t &\forall t\in \mathcal{T}\label{eq:full_model_25}\\
&\Delta{p}_t \leq |P| - 1&\forall t\in \mathcal{T}\label{eq:full_model_26}
\end{align}

\noindent \textbf{Variable definition.}
\begin{align}
&a^e_t, \Delta w^e_t, \Delta w^p_t, \Delta{p}_t\in Z^+_0 &\forall t\in \mathcal{T}\label{eq:full_model_27}\\
&d^e_t, a^p_t, d^p_t \geq 0 &\forall t \in T\label{eq:full_model_28}\\
&\alpha^r_t, \beta^r_t, \gamma^r_t  \in \{0, 1\} &\forall t\in \mathcal{T}, r\in \mathcal{R}\label{eq:full_model_29}\\
&\delta^+_t, \delta^-_t \in\{0, 1\} &\forall t\in \mathcal{T}\label{eq:full_model_30}
\end{align}

\end{document}